\documentclass[journal,10pt]{IEEEtran}     

\IEEEoverridecommandlockouts                              

\usepackage{latexsym,amssymb,amsmath, amsbsy,amsopn, amstext,amsthm}
\usepackage{graphicx,epsfig,tikz,pgf,hyperref,cancel,color,float,ifpdf}
\usepackage{extarrows,mathrsfs}
\usepackage{amsmath,  amssymb,stfloats}
\usepackage{subfigure,setspace}

\graphicspath{{./}}
\newcommand{\R}{\mathbb{R}}

\newcommand*\diff{\mathop{}\!\mathrm{d}}

\newtheorem{definition}{\bfseries Definition}
\newtheorem{proposition}{\bfseries Proposition}
\newtheorem{assumption}{\it Assumption}
\newtheorem{theorem}{\bfseries Theorem}

\newtheorem{remark}{\bfseries Remark}

\allowdisplaybreaks[4]
\title{Observer-Based Controllers for  Incrementally Quadratic Nonlinear Systems with  Disturbances}

\author{Xiangru Xu,  Beh{\c{c}}et A{\c{c}}{\i}kme{\c{s}}e, Martin J. Corless
\thanks{X. Xu is with the Department of Mechanical Engineering,  University of Wisconsin-Madison, Madison, WI, USA (email: {\tt\small xiangru.xu@wisc.edu}). }
\thanks{B. A{\c{c}}{\i}kme{\c{s}}e is with the Department of Aeronautics \& Astronautics, University of Washington, Seattle, WA, USA (email: {\tt\small behcet@uw.edu}).} 
\thanks{M. J. Corless is with the School of Aeronautics \& Astronautics, Purdue University, W. Lafayette, IN, USA (email: {\tt\small corless@purdue.edu}).}
}

\begin{document}

\maketitle
\thispagestyle{empty}
\pagestyle{empty}

\begin{abstract}
Robust global stabilization of nonlinear systems by observer-based feedback controllers is a challenging task.
This paper investigates the problem of designing observer-based stabilizing controllers for incrementally quadratic nonlinear systems with external disturbances. The nonlinearities considered in the system model satisfy the incremental quadratic constraints, which are characterized by incremental multiplier matrices and encompass many common nonlinearities. The simultaneous search for the observer and the controller gain matrices is formulated as a feasibility problem of linear matrix inequalities, for two parameterizations (i.e., the block diagonal parameterization and the block anti-triangular parameterization) of the incremental multiplier matrices, respectively. The closed-loop system implementing the observer-based feedback controller is proven to be input-to-state stable  with respect to external disturbances.  
Using the proposed continuous-time observer-based controllers, event-triggered
controllers with time regularization are constructed for globally Lipschitz systems, such that the closed-loop system is Zeno-free and input-to-state practically stable.
\end{abstract}
%

\section{Introduction}\label{sec:intro}

As the state variables of a system are difficult or expensive to measure in practice, output feedback control design has received a lot of  attention (see, e.g., \cite{krener1983linearization,vidyasagar1980stabilization,krstic1995nonlinear,andrieu2009unifying}) and found applications in biological systems \cite{borri2017luenberger,cheah2014observer,mendez2007selected}, mechanical systems \cite{khalil2017high,rajamani2017observers,su2007simple}, power systems \cite{ouassaid2012observer,mahmud2012full}, and networked control systems \cite{hong2008distributed,ahmed2012high}, among others.   
For linear systems, the output feedback stabilizing control design problem can be solved by designing the state-feedback controller and the state observer independently, which is known as the \emph{controller-observer separation principle}. For nonlinear systems, a certainty-equivalence implementation of a globally stabilizing state-feedback controller with an asymptotic observer can lead to finite escape time (e.g., see the counter-examples in \cite{kokotovic1992joy,mazenc1994global}), which makes observer-based stabilizing controller design
a challenging problem  \cite{kokotovic2001constructive}. By using a high-gain observer \cite{esfandiari1992output,khalil2014high}, separation principles for input-output linearizable systems were studied in \cite{teel1994global,khalil1993semiglobal,atassi2000separation,atassi2001separation} and semiglobal asymptotic stability of the resulting closed-loop systems was proven in these papers. Separation principles for some other special class of nonlinear systems were also investigated, such as bilinear systems \cite{gauthier1992separation,tsinias1993sontag}, non-affine nonlinear systems \cite{lin1995bounded}, systems with  nondecreasing or slope-restricted  nonlinearities \cite{arcak2001nonlinear,arcak2001observer,arcak2005certainty}, and cascaded systems \cite{fossen1999passive,loria2000separation}. 
Apart from the certainty-equivalence approach, interdependent design of the	controller and the observer was investigated in \cite{praly1993stabilization,pomet1993dynamic,praly1993lyapunov}.

Linear matrix inequalities (LMIs) provide a computationally efficient approach for the synthesis of  observer-based output feedback controllers \cite{boyd1994linear}, where the main difficulty lies in the coupling between the unknown matrices of the observer and the controller and the Lyapunov matrices. For linear systems, LMI-based conditions were proposed for the robust observer-based stabilization of linear systems with state perturbations  \cite{lien2004robust} or with parametric uncertainties \cite{kheloufi2013lmi,zemouche2017robust}. For nonlinear systems, the synthesis problem is often formulated as the feasibility of bilinear matrix inequalities (BMIs), which is known to be an NP-hard problem \cite{toker1995np}. Different approaches that aim to transform the non-convex BMI conditions to convex LMI conditions have been proposed: \cite{zemouche2017robust} studied  observer-based controller design for Lipschitz nonlinear systems with uncertain parameters, and developed an LMI-based design technique that relies on the linearization of the corresponding BMIs;  \cite{wang2018sequential} investigated  observer-based control design for the interconnection of a linear system and an uncertain nonlinear operator satisfying the integral quadratic constraint, and proposed a  sequential LMI algorithm to solve BMIs; \cite{kim2017robust} studied output feedback control
of discrete-time parametric uncertain Lure systems, and developed an LMI-based iterative algorithm to solve BMIs. Moreover, \cite{grandvallet2013new} investigated $H_\infty$ stabilization of discrete-time globally Lipschitz nonlinear systems, and provided LMI-based conditions that compute simultaneously the observer and controller gains; \cite{ekram2017observer} considered asymptotic stabilization of continuous-time Lipschitz nonlinear systems and developed LMI-based conditions that synthesize the gain matrices of the observer-based controller.

This paper considers observer-based output feedback global stabilization of a class of nonlinear systems whose nonlinearities satisfy incremental quadratic constraints. 
The incremental quadratic constraint is characterized by an incremental quadratic inequality with incremental multiplier matrices \cite{megretski1997system,acikmesethesis,acikmese2005observers,accikmecse2008stability,accikmecse2011observers,alto13incremental}. This characterization with incremental multiplier matrices provides a general framework to represent many common classes of nonlinearities (e.g., globally Lipschitz nonlinearities, incrementally sector bounded nonlinearities,  non-decreasing nonlinearities and the polytopic Jacobian nonlinearities), implying a wide range of applicability for the proposed theoretical results. 
Observer design for systems with nonlinearities satisfying incremental quadratic constraints was studied in \cite{accikmecse2011observers}, which was later generalized to the systems with  bounded exogenous disturbances in \cite{ankush2017state}. Observer-based control design for some special classes of incrementally quadratic nonlinear systems have been investigated in  \cite{arcak2001observer,ekram2017observer,chen2007robust,arcak1999observer,fan2003observer}.

Motivated by the development of networked control systems, event-triggered control (ETC) has  recently received  a lot of  attention as it provides a new control paradigm to  reduce the resource consumption of networked control systems whose communication bandwidth and computational power are usually limited  \cite{heemels2012introduction,tabuada2007event,postoyan2015framework}.  
Most of the ETC results  assume that  full-state information is available, 
but this assumption is restrictive since many systems  only have  information on their measured outputs. 
Extending results on observer-based, event-triggered control design from linear systems (e.g., see \cite{tallapragada2012event,tarbouriech2016observer,donkers2012output}) to nonlinear systems is difficult \cite{borgers2014event}.  
Existing results on ETC design for nonlinear systems mostly assume that the continuous-time observer-based controllers are already given, but the observers and controllers themselves can be hard to construct. 
When external disturbances or measurement noise are present, the triggering rules also need to be carefully  designed to rule out the Zeno phenomenon \cite{borgers2014event,dolk2017output,abdelrahim2017robust}, e.g., using time regularization to enforce a built-in lower bound for inter-execution times.

The main contributions of the paper are summarized as follows. 
For incrementally quadratic nonlinear systems affected by external disturbances and measurement noise, LMI-based sufficient conditions are developed for the design of robust stabilizing observer-based controllers. The simultaneous search for the observer and the controller gain matrices is formulated as a feasibility problem of LMIs when the incremental multiplier matrices are parameterized as the block
diagonal matrices or the block anti-triangular matrices. The resulting closed-loop system is proven to be input-to-state stable with respect to disturbances. Using the  proposed continuous-time observer-based controller, event-triggered
controllers are constructed for globally Lipschitz systems affected by external disturbances and measurement noise where the triggering rule is designed with an enforced positive lower-bound on  inter-execution times. The resulting closed-loop system is Zeno-free and input-to-state practically stable with respect to external disturbances.	
A preliminary version of this work appeared in  \cite{xu2018acc}. The present paper is different from 
\cite{xu2018acc} in the following important ways:  the system model considered is subject to external disturbances and measurement noise; the event-triggered mechanism is considered;  complete proofs are included, and more discussion is added.  
The remainder of the paper is organized as follows: Section \ref{sec:formulation} introduces preliminaries on incremental quadratic constraints and input-to-state practical stability,  Section \ref{sec:continuous} develops LMI-based conditions for the design of robust stabilizing observer-based controllers for two parameterizations of the incremental multiplier matrices, 
Section \ref{sec:trigger} presents the event-triggered controller design, 
Section \ref{sec:example} provides a simulation example, and Section \ref{sec:conclusion} provides the conclusions.

\emph{Notation.} $\R_0^+$ denotes the set of non-negative real numbers; $\|x\|$ denotes the $2$-norm of a vector $x$;  
$\|P\|$ denotes the maximum singular value of a matrix $P$; 
$\lambda_{m}(P)$ and $\lambda_{M}(P)$ denote the minimum and maximum eigenvalues of a symmmetric matrix $P$, respectively; $I_n$ denotes an identity matrix of size $n$; ${\bf 0}_{n_1\times n_2}$ and ${\bf 0}_{n}$ denote the zero matrix of size  $n_1\times n_2$ and 
the zero vector of size $n$, respectively, 
where the subscript will be omitted when  clear from context. 
For symmetric matrices, $*$ denotes entries whose values follow from symmetry.
For a matrix $M$, $M\succ 0$, $M\succeq 0$, $M\prec 0$, $M\preceq 0$ mean $M$ is positive definite,   positive semi-definite, negative definite, and negative semi-definite,  respectively. 
A continuous function $f: \R_0^+\rightarrow \R_0^+$ belongs to class $\mathcal{K}$ (denoted as $f\in\mathcal{K}$) if it is strictly increasing and $f(0)=0$; $f$ belongs to class $\mathcal{K}_{\infty}$ (denoted as $f\in\mathcal{K}_{\infty}$) if $f\in\mathcal{K}$ and $f(r)\rightarrow\infty$ as $r\rightarrow \infty$. A continuous function $f: \R^+_0\times \R^+_0\rightarrow \R^+_0$ belongs to class $\mathcal{KL}$ (denoted as $f\in\mathcal{KL}$) if for each fixed $s$, the function $f(\cdot,s)\in \mathcal{K}_{\infty}$  and for each each fixed $r$, the function $f(r,\cdot)$ is decreasing  and $f(r,s)\rightarrow 0$ as $s\rightarrow 0$.

\section{Preliminaries}\label{sec:formulation}

Consider the following nonlinear system
\begin{align}\label{dyn1}
\begin{cases}
\dot x=Ax+Bu+Ep(q)+E_ww,\\
y=Cx+Du+F_ww,\\
q=C_qx,
\end{cases}
\end{align}
where $x\in\R^{n_x}$ is the state, $u\in\R^{n_u}$ is the control input, $y\in\R^{n_y}$ is the measured output, $p:\R^{n_q}\rightarrow \R^{n_p}$ is the known nonlinearity of the system, $w\in\R^{n_w}$ is the unknown external disturbance or measurement noise, and $A\in\R^{n_x\times n_x},B\in\R^{n_x\times n_u},C\in\R^{n_y\times n_x},D\in\R^{n_y\times n_u},C_q\in\R^{n_q\times n_x},E\in\R^{n_x\times n_p},E_w\in\R^{n_x\times n_w},F_w\in\R^{n_y\times n_w}$ are constant matrices of appropriate dimensions. 

The characterization of the nonlinearity $p$ is based on incremental multiplier matrices \cite{accikmecse2011observers,alto13incremental}.
\begin{definition}\label{def:delQC}
	Given a function $p:\R^{n_q}\rightarrow \R^{n_p}$, a symmetric matrix $M\in\R^{(n_q+n_p)\times (n_q+n_p)}$ is called an \emph{incremental multiplier matrix} ($\delta$-MM)
	for $p$ if it satisfies the following \emph{incremental quadratic constraint} for all $q_1, q_2\in\R^{n_q}$:
	\begin{equation}\label{eq:delQC}
	\begin{pmatrix}
	\delta q\\
	\delta p
	\end{pmatrix}^\top M 
	\begin{pmatrix}
	\delta q\\
	\delta p
	\end{pmatrix}\geq  0
	\end{equation} 
	where $\delta q=q_2-q_1$, $\delta p=p(q_2)-p(q_1)$.	
\end{definition}

	For a given nonlinearity $p$, its $\delta$-MM is not unique. Denote $\mathcal{M}$ as the set of incremental multiplier matrices for $p$. If $M\in\mathcal{M}$, then  $\lambda M\in\mathcal{M}$ for any $\lambda \ge 0$.

	\begin{remark}\label{remark1}
The global Lipschitz condition  $\|p(q_2)-p(q_1)\|\leq \gamma \|q_2-q_1\|$
		where $\gamma\!>\!0$  can be expressed in the form of \eqref{eq:delQC} with  
		\begin{align}\label{M1}
		M= 
		\begin{pmatrix}
		\gamma^2 I&{\bf 0}\\
		{\bf 0}&-I
		\end{pmatrix}.
		\end{align} 
		The incrementally sector bounded nonlinearity $(\delta p-K_1\delta q)^\top S(\delta p-K_2\delta q)\leq 0$ 
		where $S=S^\top$ can be expressed in the form of \eqref{eq:delQC} with 
		\begin{align}\label{M2}
		M=
		\begin{pmatrix}
		-K_1^\top SK_2-K_2^\top SK_1&*\\
		S(K_1+K_2)&-2S
		\end{pmatrix}.
		\end{align} 
		The nondecreasing nonlinearity, which satisfies $\delta p^\top \delta q\geq 0$,
		can be expressed in the form of \eqref{eq:delQC} with 
		\begin{align}\label{M3}
		M=
		\begin{pmatrix}
		{\bf 0}&I\\
		I&{\bf 0}
		\end{pmatrix}.
		\end{align} 
		Refer to \cite{accikmecse2011observers,alto13incremental} for some other nonlinearities that can be expressed using the incremental quadratic constraint.
\end{remark}

Next, we introduce input-to-state practical stability and its characterization using Lyapunov functions. Consider the system  
\begin{align}\label{eqnISPS}
\dot{x}=f(x,u)
\end{align}
where $f:\R^{n_x}\times \R^{n_u}\rightarrow \R^{n_x}$ is a locally Lipschitz function and $u:\R\rightarrow \R^{n_u}$ is a measurable essentially bounded input. Define $x(t,x_0,u)$ as the solution of \eqref{eqnISPS} with initial state $x_0$ and input $u$, which satisfies $x(0,x_0,u)=x_0$.

\begin{definition}\label{dfn:ISPS}(Def. 2.1 of \cite{jiang1996lyapunov})
	The system \eqref{eqnISPS} is called \emph{input-to-state practically stable} (ISpS) w.r.t. $u$, if there exist functions $\beta_1\in\mathcal{KL}$, $\beta_2\in\mathcal{K}$ and a non-negative constant $d$ such that for every initial state $x_0$ and every measurable essentially bounded $u$ defined  on $[0,\infty)$, the solution $x(t,x_0,u)$ exists on $[0,\infty)$ and satisfies 
	\begin{align}\label{ISSineq3}
	\|x(t,x_0,u)\|\leq \beta_1(\|x_0\|,t)+\beta_2(\|u\|_{\infty})+d,\;\forall t\geq 0
	\end{align}
	where $\|u\|_{\infty}:=\emph{ess sup}_{t\geq 0}\|u(t)\|$.
\end{definition}
When \eqref{ISSineq3} is satisfied with $d=0$, the system is said to be \emph{input-to-state stable} (ISS) w.r.t. $u$~\cite{isidori2013nonlinear}. 

\begin{definition}(Remark 2.2 of \cite{jiang1996lyapunov})
	A smooth function $V:\R^n\rightarrow \R$ is said to be an \emph{ISpS-Lyapunov function} for the system \eqref{eqnISPS} if   $V$ is  radially unbounded, positive definite  and there exist functions $\gamma\in\mathcal{K}_\infty,\chi\in\mathcal{K}$ and a non-negative constant $d$ such that the following condition holds:
	\begin{align}\label{ISSineq}
	\nabla V(x)^Tf(x,u)\leq -\gamma(\|x\|)+\chi(\|u\|)+d.
	\end{align} 
\end{definition}

Instead of requiring inequality \eqref{ISSineq}, the ISpS-Lyapunov function can be also defined equivalently as follows: a smooth, positive definite, radially unbounded function $V$ is an ISpS-Lyapunov function for the system \eqref{eqnISPS} if there exist a positive-definite function $\gamma$, a class $\mathcal{K}$ function $\chi$ and a non-negative constant $d$ such that the following condition holds (Def. 2.2 of \cite{jiang1996lyapunov}):
\begin{align}
\|x\|\geq \chi(\|u\|)+d\; \Rightarrow \; \nabla V(x)^Tf(x,u)\leq -\gamma(\|x\|).\label{ISSineq4}
\end{align}

The existence of an \emph{ISpS-Lyapunov function} is a necessary and sufficient condition for the ISpS property.

\begin{proposition}\label{proISS}\cite{sontag1995characterizations}
	The system \eqref{eqnISPS} is ISpS (resp. ISS) if and only if it has an ISpS- (resp. ISS-) Lyapunov function.
\end{proposition}

In particular, if there exist a symmetric and positive definite matrix $P=P^\top\succ0$,  two constants $\alpha>0,d\geq0$ and a function $\chi\in\mathcal{K}_\infty$ such that the positive definite function $V(x)=x^\top Px$ satisfies
\begin{align}\label{ISSineq2}
\nabla V(x)^Tf(x,u)\leq -\alpha V(x)+\chi(\|u\|)+d,
\end{align} 
then $V$ is an ISpS-Lyapunov function satisfying \eqref{ISSineq} with $\gamma(\|x\|)=\alpha\lambda_{m}(P)\|x\|^2$, implying that \eqref{eqnISPS} is ISpS w.r.t. $u$.

\section{LMI-based Conditions For Robust Global Stabilization of Incrementally Quadratic Nonlinear Systems}\label{sec:continuous}

Consider a system described by \eqref{dyn1} where the nonlinear  term $p$ satisfies  the incremental quadratic constraint \eqref{eq:delQC} for some $M\in\mathcal{M}$. In this section, a continuous-time observer and a feedback controller will be designed for \eqref{dyn1}, such that the closed-loop system is ISS w.r.t. $w$.  LMI-based sufficient conditions will be given for the simultaneous design of the observer and controller gain matrices.

The  following observer is proposed:
\begin{align}\label{obser1}
\begin{cases}
\dot{\hat{x}}\!=\!A\hat x+Bu+E p(\hat q+L_1(\hat y-y))+L_2(\hat y-y),\\
\hat y\!=\!C\hat x+Du,\\
\hat q\!=\!C_q\hat x,
\end{cases}
\end{align}
where $L_1,L_2$ are gain matrices to be designed. 
This observer contains a copy of the plant and two correction terms, the  \emph{nonlinear injection term} $L_1(\hat y-y)$ and the Luenberger-type \emph{correction term} $L_2(\hat y-y)$.
Based on observer \eqref{obser1}, we design the feedback controller $u$ as 
\begin{align}\label{input1}
u&=k(\hat x)
\end{align}
where $k:\R^{n_x}\rightarrow \R^{n_u}$ is a function that 
has the  form of
\begin{align}\label{inputform}
k(x)&=K_1x+K_2p(C_qx)
\end{align}  
with gain matrices $K_1\in\R^{n_u\times n_x}$, $K_2\in\R^{n_u\times n_p}$ to be designed.
Defining the estimation error by 
$$
e(t)=x(t)-\hat x(t)
$$
the input \eqref{input1} can be rewritten as 
$$
u=k(x)-\Delta k(x,\hat x)
$$ 
where $\Delta k(x,\hat x)=k(x)-k(\hat{x})$. 
Recalling \eqref{inputform}, $\Delta k$ can be expressed as $\Delta k=K_1e-K_2\Delta p$ 
where 
\begin{align}
\Delta p=p(\hat q)-p(q).\label{Deltap}
\end{align}

%
The closed-loop system resulting from  the observer-based controller \eqref{input1} can now be expressed as
\begin{align}\label{dyn3-1}
\begin{cases}
\dot x\!=\!(A\!+\!BK_1)x\!+\!(E\!+\!BK_2)p\!-\!B\Delta k\!+\!E_ww, \\
\dot e\!=\!(A\!+\!L_2C)e\!-\!E\delta p\!+\!(E_w\!+\!L_2F_w)w, 
\end{cases}
\end{align}
where
\begin{align}\label{delp}
\begin{cases}
\delta p=p(q+\delta q)-p(q),\\
\delta q=-(C_q+L_1C)e-L_1F_ww.
\end{cases}
\end{align}
Defining $z=
\begin{pmatrix}
x\\
e
\end{pmatrix}$, dynamics \eqref{dyn3-1} are expressed compactly as
\begin{align}
\dot z=A_cz+H_1p+H_2\delta p+H_3\Delta p +H_4w\label{eqthm3-1}
\end{align}
where $\Delta p$ is given in \eqref{Deltap}, $\delta p$ is given in \eqref{delp}, and
\begin{align}
A_c&=
\begin{pmatrix}
A+BK_1&-BK_1\\
{\bf 0}& A+L_2 C
\end{pmatrix},\label{eqA}
\end{align}
\begin{align}\label{H1}
\begin{cases}
H_1=
\begin{pmatrix}
E+BK_2\\
{\bf 0}
\end{pmatrix},\;H_2=
\begin{pmatrix}
{\bf 0}\\
-E
\end{pmatrix},\\
H_3=
\begin{pmatrix}
BK_2\\
{\bf 0}
\end{pmatrix},\;H_4=
\begin{pmatrix}
E_w\\
E_w+L_2F_w
\end{pmatrix}.
\end{cases}
\end{align}	

The following proposition provides a sufficient condition for the closed-loop system \eqref{eqthm3-1} to be ISS w.r.t. $w$.

\begin{proposition}\label{prp1}
	Consider the system described by \eqref{dyn1}-\eqref{eq:delQC} with $p({\bf 0})={\bf 0}$. Suppose that 
	there exist matrices $L_1\in\R^{n_q\times n_y},L_2\in\R^{n_x\times n_y},K_1\in\R^{n_u\times n_x},K_2\in\R^{n_u\times n_p},P\in\R^{2n_x\times 2n_x}$ with $P\succ 0$, and real numbers $\alpha_0>0,\mu>0,\sigma_1\geq 0,\sigma_2\geq 0,\sigma_3\geq 0$ such that
	\begin{align}
	\begin{pmatrix}
	S_0&S_1\\
	\ast&{\bf 0}
	\end{pmatrix}
	&+\sigma_1S_2^{\top} MS_2+\sigma_2  S_3^{\top} MS_3
	+\sigma_3S_4^{\top}MS_4 \nonumber\\ 
	& -\mu S_5^\top S_5\preceq 0\label{recomLMI}
	\end{align}
	where
	\begin{align*}
	S_0&=PA_c+A_c^\top P+\alpha_0 P,        \\
	S_1&=\left(PH_1\;\; PH_2\;\;PH_3\;\; PH_4\right),             \\
	S_2&=
	\begin{pmatrix}
	C_q \qquad 
	{\bf 0}_{n_q\times (n_x+3n_p+n_w)}\\
	{\bf 0}_{n_p\times 2n_x} \quad I_{n_p}
	\quad {\bf 0}_{n_p\times (2n_p+n_w)}
	\end{pmatrix},
	\\
	S_3&=
	\begin{pmatrix}
	{\bf 0}_{n_q\times n_x}\quad  -(C_q+L_1C)
	\quad {\bf 0}_{n_q\times 3n_p}      \quad -L_1F_w\\
	{\bf 0}_{n_p\times (2n_x+n_p)} \quad  I_{n_p}
	\quad {\bf 0}_{n_p\times (n_p+n_w)}
	\end{pmatrix},
	\\
	S_4&=
	\begin{pmatrix}
	{\bf 0}_{n_q\times n_x}\quad -C_q \quad {\bf 0}_{n_q\times (3n_p+n_w)}\\
	{\bf 0}_{n_p\times (2n_x+2n_p)}\quad I_{n_p}\quad {\bf 0}_{n_p\times n_w}
	\end{pmatrix},
	\\
	S_5&=({\bf 0}_{n_w\times (2n_x+3n_p) } \quad I_{n_w}).
	\end{align*}
	Then the closed-loop system \eqref{eqthm3-1} is ISS w.r.t. $w$ and satisfies
		$\dot V\leq -\alpha_0 V+\mu \| w\|^2$
		where $V(z)=z^\top Pz$.
	
\end{proposition}

\begin{proof}
	Since $M$ is a $\delta$-MM for $p$ and $p({\bf 0})={\bf 0}$, 
	it holds that
	$$
	\begin{pmatrix}
	q\\
	p
	\end{pmatrix}^\top\!\!M 
	\begin{pmatrix}
	q\\
	p
	\end{pmatrix}\!\geq\! 0,
	\begin{pmatrix}
	\delta q\\
	\delta p
	\end{pmatrix}^\top\!\!M 
	\begin{pmatrix}
	\delta q\\
	\delta p
	\end{pmatrix}
	\!\geq\! 0,
	\begin{pmatrix}
	\Delta q\\
	\Delta p
	\end{pmatrix}^\top\!\!M 
	\begin{pmatrix}
	\Delta q\\
	\Delta p
	\end{pmatrix}\!\geq\! 0,
	$$
	where $\delta p,\delta q$ are given in \eqref{delp}, $\Delta p$ is given in \eqref{Deltap}, and 
	\begin{align}\label{Deltaq}
	\Delta q=C_q\hat x-C_qx=-C_qe.
	\end{align}
	With  $\xi=(x^\top\: e^\top\:  p^\top\: \delta p^\top\: \Delta p^\top\: w^\top)^\top$,
	$$
	\begin{pmatrix}
	q\\
	p
	\end{pmatrix}=S_2\xi, \quad
	\begin{pmatrix}
	\delta q\\
	\delta p
	\end{pmatrix}=S_3\xi,\quad
	\begin{pmatrix}
	\Delta q\\
	\Delta p
	\end{pmatrix}=S_4\xi.$$
	Hence, 
	$\xi^\top S_2^\top MS_2 \xi\geq 0,\xi^\top S_3^\top MS_3 \xi \geq 0 ,\xi^\top S_4^\top MS_4 \xi\geq 0$. 
	Pre- and post-multiply \eqref{recomLMI} by $\xi^\top$ and $\xi$, respectively. Since
	$\sigma_1,\sigma_2,\sigma_3$ are non-negative, we obtain that
	\begin{align}
	\xi^\top
	\begin{pmatrix}
	S_0&S_1\\
	\ast &{\bf 0}
	\end{pmatrix}\xi-\mu \xi^\top S_5^\top S_5\xi\leq 0.\label{ineqpro1}
	\end{align}
	Consider the positive definite function defined by $V(z)=z^\top Pz$. Then, it is easy to check that $\dot V+\alpha_0 V-\mu w^\top w$ is equal to the left hand side of \eqref{ineqpro1} where $\dot V$ is the derivative of $V$ along the trajectories of \eqref{eqthm3-1}. Therefore, $V$ is an ISS-Lyapunov function since $\dot V\leq -\alpha_0 V+\mu \| w\|^2$. The conclusion follows from Proposition \ref{proISS}.	
\end{proof}
Clearly,  matrix inequality \eqref{recomLMI} is not a LMI. 
In the next two subsections, we will consider two parameterizations of the $\delta$-MM $M$ and provide  LMI conditions which can be used to solve for $M$ and gain matrices $L_1,L_2,K_1,K_2$ simultaneously.

\subsection{Block Diagonal Parameterization}

This subsection considers a block diagonal parameterization of the $\delta$-MM for $p$. 
We first make the following two assumptions on the  parameterizations of $M$.

	\begin{assumption}\label{ass1}
		There exist  a set $\mathcal{N}_1$ of matrix pairs $(X_1,Y_1)$ with $X_1\in\R^{n_q\times n_q}$, $Y_1\in\R^{n_p\times n_p}$ symmetric,  and an invertible matrix $T_1$ with
		\begin{align}\label{eqT}
		T_1=
		\begin{pmatrix}
		T_{11}&T_{12}\\
		T_{13}&T_{14}
		\end{pmatrix}
		\end{align}
		and $T_{14} \in\R^{n_p\times n_p}$  invertible, such that $M_1$ given below is a $\delta$-MM of $p$  for all  $(X_1,Y_1)\in\mathcal{N}_1$:
		\begin{align}
		\label{eq:Mass1}
		M_1&=T_1^\top \tilde M_1T_1\;\;\mbox{where}\;\;\tilde M_1=
		\begin{pmatrix}
		X_1&{\bf 0}\\
		{\bf 0}&-Y_1
		\end{pmatrix}.
		\end{align}  
	\end{assumption}

	\begin{assumption}\label{ass2}
		There exist  a set $\mathcal{N}_2$ of matrix pairs 
		$(X_2,Y_2)$ with $X_2\in\R^{n_q\times n_q}$, $Y_2\in\R^{n_p\times n_p}$ symmetric and invertible,  and an invertible matrix $T_2$ with
		\begin{align}\label{eqT2}
		T_2=
		\begin{pmatrix}
		T_{21}&T_{22}\\
		T_{23}&T_{24}
		\end{pmatrix}
		\end{align}
		and $T_{24} \in\R^{n_p\times n_p}$  invertible, such that $M_2$ given below is a $\delta$-MM of $p$ for all $(X_2,Y_2)\in\mathcal{N}_2$:
		\begin{align}
		\label{eq:Mass2}
		M_2&=T_2^\top \tilde M_2T_2\;\;\mbox{where}\;\;\tilde M_2=
		\begin{pmatrix}
		X_2^{-1}&{\bf 0}\\
		{\bf 0}&-Y_2^{-1}
		\end{pmatrix}.
		\end{align}  
	\end{assumption}

\begin{remark}\label{remark2}
	For the globally Lipschitz nonlinearity $\|p(q_2)-p(q_1)\|\leq \gamma \|q_2-q_1\|$, the matrix $M$ in \eqref{M1} satisfies Assumption \ref{ass1} and \ref{ass2} if we choose
	\begin{align*}
	T_1\!=\!T_2\!=\!
	\begin{pmatrix}
	\gamma I&{\bf 0}\\
	{\bf 0}&I
	\end{pmatrix},\;\mathcal{N}_1\!=\!\mathcal{N}_2\!=\!\{(\lambda I,\lambda I)|\lambda>0\}.
	\end{align*}
	For the incrementally sector bounded nonlinearity $(\delta p-K_1\delta q)^\top S(\delta p-K_2\delta q)\leq 0$ where $S$ is symmetric and invertible, the  matrix $M$ in \eqref{M2}  satisfies Assumption \ref{ass1} and \ref{ass2} if we choose 
	\begin{align*}
	T_1\!=\!T_2\!=\!
	\begin{pmatrix}
	K_2-K_1&{\bf 0}\\
	K_2+K_1&-2I
	\end{pmatrix},\;
	\mathcal{N}_1\!=\!\mathcal{N}_2\!=\!\{(\lambda S,\lambda S)|\lambda>0\}.
	\end{align*}
	For the nondecreasing nonlinearity $\delta p^\top \delta q\geq 0$, the matrix $M$ in \eqref{M3} satisfies Assumption \ref{ass1} and \ref{ass2} if we choose
	\begin{align*}
	T_1\!=\!T_2\!=\!
	\begin{pmatrix}
	I&I\\
	I&-I
	\end{pmatrix},\;
	\mathcal{N}_1\!=\!\mathcal{N}_2\!=\!\{(\lambda I,\lambda I)|\lambda>0\}.
	\end{align*}
		$\mathcal{N}_1$ and $\mathcal{N}_2$ do not have to be the set of scalings of a matrix pair as in the examples above. For instance, for the nonlinearity whose Jacobian is confined within a polytope or a cone, $\mathcal{N}_1$ that satisfies Assumption \ref{ass1} (or $\mathcal{N}_2$ that satisfies Assumption \ref{ass2}) is characterized  via matrix inequalities (see Section 5 in \cite{accikmecse2011observers} for more details). Furthermore, $T_1$ does not necessarily has to be chosen to be equal to  $T_2$.
\end{remark}

Because $T_1$ in Assumption \ref{ass1} and $T_2$ in Assumption \ref{ass2} are invertible, the matrix $\Gamma_{i1} (i=1,2)$ defined as
\begin{align}
\Gamma_{i1}&=T_{i1}-T_{i2}T_{i4}^{-1}T_{13}\label{gamma1}
\end{align}
is also invertible by the matrix inversion lemma. Furthermore, we define the matrix $\Gamma_{i2} (i=1,2)$ as
\begin{align}
\Gamma_{i2}=T_{i2}T_{i4}^{-1}.\label{gamma2}
\end{align}

The following theorem  provides sufficient conditions for the design of matrices $L_1,L_2$ in the observer \eqref{obser1} and matrices $K_1,K_2$ in the controller \eqref{input1}, when the $\delta$-MM can be parameterized in a block diagonal manner.

\begin{theorem}\label{thm1}
	Consider the system described by \eqref{dyn1}-\eqref{eq:delQC}  with  $p({\bf 0})={\bf 0}$. Suppose that\\
		1)
		Assumption \ref{ass1} holds; \\
		2)
		Assumption \ref{ass2} holds with $M_2=
		\begin{pmatrix}
		M_{21}&M_{22}\\
		M_{23}&M_{24}
		\end{pmatrix}$
		where $M_{24}\in\R^{n_p\times n_p}$ and $M_{24}\prec 0$; \\
		3) 
		there exist positive numbers $\alpha_1,\alpha_2,\mu_1,\mu_2$,  matrices $R_1,R_2,R_3,R_4$, symmetric  and positive definite matrices
		$P_1, P_2, X_1, X_2,Y_2$ and  a symmetric matrix $Y_1$,  such that $(X_1,Y_1)\in\mathcal{N}_1$, $(X_2,Y_2)\in\mathcal{N}_2$ and
	\begin{align}
	\text{\rm (observer ineq.)}\;\quad  &  
	\begin{pmatrix}
	\Phi-\varphi^\top Y_1\varphi & \phi^\top\\
	\phi & -X_1
	\end{pmatrix}\preceq 0, \label{LMI1}\\
	\text{\rm (controller ineq.)} \;\quad &
	\begin{pmatrix}
	\Psi-\varphi^\top Y_2\varphi & \psi^\top\\
	\psi & -X_2
	\end{pmatrix}\preceq 0, \label{LMI2}
	\end{align} 
	where
	\begin{align}
	\Phi&=
	\begin{pmatrix}
	\Phi_0 & -P_1\tilde E_1 & P_1 E_w+R_1 \\
	* & {\bf 0} & {\bf 0} \\
	* &*&-\mu_1I
	\end{pmatrix},\label{Phi}\\
	\Phi_0&=\tilde A_1^\top P_1+P_1\tilde A_1+ C^\top R_1^\top+R_1 C+\alpha_1 P_1,\label{Phi0}\\
	\Psi&=
	\begin{pmatrix}
	\Psi_0& \tilde E_2Y_2 +BR_4 & E_w \\
	*  & {\bf 0} &{\bf 0}  \\
	* & *&-\mu_2I
	\end{pmatrix},\label{Psi}\\
	\Psi_0&=\tilde A_2 P_2+P_2\tilde A_2^\top+ B R_3+R_3^\top B^\top+\alpha_2 P_2,\label{Psi0}\\
	\phi&=(-(X_1 \Gamma_{11} C_q+R_2 C) , X_1 \Gamma_{12}, -R_2F_w),\label{phi}\\
	\varphi&=({\bf 0}_{n_p\times n_x}, I_{n_p}, {\bf 0}_{n_p\times n_w}),\label{varphi}\\
	\psi&=(\Gamma_{21} C_q P_2 , \Gamma_{22}Y_2 ,{\bf 0}_{n_q\times n_w}),\\
	\tilde A_i&=A-ET_{i4}^{-1}T_{i3}C_q, \;i=1,2,\label{tilA}\\
	\tilde E_i&=ET_{i4}^{-1},\;i=1,2, \label{tilE}
	\end{align}
	with  $\Gamma_{i1}(i=1,2)$ given in \eqref{gamma1} and $\Gamma_{i2}(i=1,2)$ given in \eqref{gamma2}.
	Then, the  closed-loop system \eqref{dyn3-1} is ISS w.r.t. $w$ with
	\begin{align}\label{L1}
	\begin{cases}
	L_1=\Gamma_{11}^{-1}X_1^{-1}R_2,\\
	L_2=P_1^{-1}R_1+ET_{14}^{-1}T_{13}L_1,\\
	K_1=R_3P_2^{-1} + K_2T_{24}^{-1}T_{23}C_q,\\
	K_2=R_4 Y_2^{-1}T_{24}.
	\end{cases}
	\end{align} 
\end{theorem}
\begin{proof}
	The proof proceeds in five steps.
	
	1) Firstly, we derive dynamics of the system under  transformations of variables $q$ and $p$ via $T_1$ and $T_2$. 
	Since $M_1$ (resp. $M_2$) satisfies Assumption \ref{ass1} (resp. Assumption \ref{ass2}) with an invertible matrix $T_1$ (resp. $T_2$),
	we introduce variable transformations from $(q, p)$ to $(\tilde q_i,\tilde p_i)$ as follows:
	\begin{align}\label{tilq}
	\begin{pmatrix}
	\tilde q_i\\
	\tilde p_i
	\end{pmatrix}=T_i
	\begin{pmatrix}
	q\\
	p
	\end{pmatrix}.
	\end{align} 
	Since 
	$
	\tilde{p}_i=T_{i3}q+T_{i4}p
	$
	and $T_{i4}$ is invertible, we have $p=T_{i4}^{-1}\tilde{p}_i-T_{i4}^{-1}T_{i3}q$ 
	and $\tilde{q}_i=\Gamma_{i1}q+\Gamma_{i2}\tilde p_i$ for $i=1,2$, where $\Gamma_{i1},\Gamma_{i2}$ are given in \eqref{gamma1},\eqref{gamma2}.
	Recall that  $\Gamma_{i1}$ is invertible since $T_i$ is invertible. 
	
	Substituting $p=T_{24}^{-1}\tilde{p}_2-T_{24}^{-1}T_{23}q$ into  \eqref{dyn3-1}, we have
	\begin{align}
	\dot x&=(\tilde A_2+B\tilde{K}_1)x+(\tilde E_2+B\tilde{K}_2)\tilde p_2-B\Delta k+E_ww,\label{dyn3-2}
	\end{align}
	where $\tilde p_2=\tilde p_2(C_qx)$, $\tilde{A}_2$ is given in \eqref{tilA}, $\tilde E_2$ is given in \eqref{tilE},
	\begin{align}
	\tilde{K}_1 &= K_1 -K_2T_{24}^{-1}T_{23}C_q,\qquad \tilde{K}_2 = K_2T_{24}^{-1},\label{tildeK}
	\end{align}
	and
	\begin{equation}
	\Delta k=\tilde{K}_1e-\tilde{K}_2\Delta \tilde p,\label{Delk}
	\end{equation}
	where 
	\begin{align}
	\Delta \tilde{p}&=\tilde{p}_2 (C_q\hat{x})-\tilde{p}_2(C_qx).\label{Delp}
	\end{align}
	Define $\begin{pmatrix}\delta\tilde q_1\\
	\delta\tilde p_1
	\end{pmatrix}=T_1
	\begin{pmatrix}
	\delta q\\
	\delta p
	\end{pmatrix}$. Then, $\delta\tilde p_1=T_{13}\delta q+T_{14}\delta p$, which implies that $\delta p=T_{14}^{-1}\delta \tilde p_1-T_{14}^{-1}T_{13}\delta q$.   
	Substituting this form of $\delta p$ and \eqref{delp} into \eqref{dyn3-1}, we have
	\begin{align}
	\dot e&=(\tilde A_1+\tilde L_2C)e-\tilde E_1\delta\tilde p_1+(E_w+\tilde L_2F_w)w,\label{dyn4-2}
	\end{align}
	where $\tilde{A}_1$ is given in \eqref{tilA}, $\tilde E_1$ is given in \eqref{tilE}, and $\tilde L_2$ is defined as
	\begin{align}
	\tilde L_2&=L_2-ET_{14}^{-1}T_{13}L_1.\label{tilL2}
	\end{align}
	
	Equations \eqref{dyn3-2} and \eqref{dyn4-2} are the dynamics of the closed-loop system after transformations of variables via $T_1$ and $T_2$.

	2) We now consider the performance of the observer.
	From \eqref{L1} we have $R_1=P_1\tilde L_2$ where $\tilde L_2$ is given in \eqref{tilL2}, and  $R_2=X_1\Gamma_{11}L_1$. Plugging $R_1$ into $\Phi$ in \eqref{Phi}, we have $\Phi_0=P_1(\tilde A_1+\tilde L_2C)+(\tilde A_1+\tilde L_2C)^\top P_1+\alpha_1 P_1$, and the $(1,3)$ entry of $\Phi$ to be  $P_1 E_w+R_1F_w=P_1(E_w+\tilde L_2F_w)$; plugging $R_2$ into $\phi$ in \eqref{phi} we have $\phi=X_1\phi_0$ where $\phi_0:=(-\Gamma_{11} (C_q+L_1 C) ,  \Gamma_{12}, -\Gamma_{11}L_1F_w).$ 
	Recalling $\varphi$ given in \eqref{varphi} and applying Schur's complement to \eqref{LMI1}, we have
	\begin{align}
	&\Phi+
	\begin{pmatrix}
	\phi_0 \\
	\varphi
	\end{pmatrix}^\top\tilde M_1
	\begin{pmatrix}
	\phi_0 \\
	\varphi
	\end{pmatrix}\preceq 0.\label{ineqV1}
	\end{align}

	Define $\xi_1=(e^\top,\delta\tilde p_1^\top,w^\top)^\top$. Pre- and post-multiplying the inequality \eqref{ineqV1} by $\xi_1^\top$ and $\xi_1$, respectively, we have 
	\begin{align}
	\xi_1^\top \Phi\xi_1+\xi_1^\top
	\begin{pmatrix}
	\phi_0 \\
	\varphi
	\end{pmatrix}^\top\tilde M_1
	\begin{pmatrix}
	\phi_0 \\
	\varphi
	\end{pmatrix}\xi_1\leq 0.\label{ineqV11}
	\end{align}
	
Note that $\delta\tilde q_1=T_{11}\delta q+T_{12}\delta p=T_{11}\delta q+T_{12}T_{14}^{-1}\delta \tilde p_1-T_{12}T_{14}^{-1}T_{13}\delta q=\Gamma_{11}\delta q+\Gamma_{12}\delta\tilde p_1=-\Gamma_{11}(C_q+L_1C)e+\Gamma_{12}\delta\tilde p_1-\Gamma_{11}L_1F_ww.$ Therefore, $\begin{pmatrix}
\delta\tilde q_1\\
\delta\tilde p_1
\end{pmatrix}=\begin{pmatrix}
\phi_0 \\
\varphi
\end{pmatrix}\xi_1$. 
	Since $\begin{pmatrix}
	\delta q\\
	\delta p
	\end{pmatrix}^\top M
	\begin{pmatrix}
	\delta q\\
	\delta p
	\end{pmatrix}\geq 0$, we have $\begin{pmatrix}
	\delta\tilde q_1\\
	\delta\tilde p_1
	\end{pmatrix}^\top \tilde M_1 
	\begin{pmatrix}
	\delta\tilde q_1\\
	\delta\tilde p_1
	\end{pmatrix}\geq 0$, and therefore,
	$\xi_1^\top 
	\begin{pmatrix}
	\phi_0 \\
	\varphi
	\end{pmatrix}^\top \tilde M_1 
	\begin{pmatrix}
	\phi_0 \\
	\varphi
	\end{pmatrix}\xi_1\geq 0$. 
	Thus, $\xi_1^\top \Phi\xi_1\leq 0$ from \eqref{ineqV11}, which is equivalent to $2e^\top P_1[(\tilde A+\tilde L_2\tilde C)e-\tilde E_1 \delta\tilde p_1+(E_w+\tilde L_2F_w)w]
	+\alpha_1 e^\top P_1e-\mu_1\|w\|^2\leq 0.$

	Define $V_1(e)=e^\top P_1e$. Then the derivative of $V_1$ along the trajectory of \eqref{dyn4-2} satisfies 
	\begin{align}
	\dot V_1&=2e^\top P_1[(\tilde A+\tilde L_2\tilde C)e-\tilde E_1 \delta\tilde p_1+(E_w+\tilde L_2F_w)w]\nonumber\\
	&\leq -\alpha_1 e^\top P_1e+\mu_1\|w\|^2.\label{eqethm1}
	\end{align}
	
	3) We now prove that $\|\Delta k\|/\|e\|$ is bounded where $\Delta k$ is given in \eqref{Delk}.
	Since $M_{24}=T_{22}^\top X_2^{-1}T_{22}-T_{24}^\top Y_2^{-1}T_{24}\prec0$ and $T_{24}$ is invertible, we have
	\begin{align}
	\Gamma_{22}^\top X_2^{-1} \Gamma_{22}-Y_2^{-1}=T_{24}^{-\top}M_{24}T_{24}^{-1}\prec 0.\label{Mpf1}
	\end{align}
	
	Recall that $\Delta q=-C_qe$ in \eqref{Deltaq} and define 
	$\Delta \tilde q:=\tilde q_2(C_q\hat x)-\tilde q_2(C_qx)$. 
	Then, $\Delta \tilde q=-\Gamma_{21}C_qe+\Gamma_{22}\Delta\tilde p$ where $\Delta\tilde p$ is given in \eqref{Delp}. Define $\zeta=(e^\top,\Delta\tilde p^\top)^\top$. Therefore, 
	\begin{align}
	&\zeta^\top
	\begin{pmatrix}
	-\Gamma_{21}C_q & \Gamma_{22} \\
	{\bf 0} & I
	\end{pmatrix}^\top\tilde M_2
	\begin{pmatrix}
	-\Gamma_{21}C_q & \Gamma_{22} \\
	{\bf 0} & I
	\end{pmatrix}\zeta\nonumber\\
	=& 
	\begin{pmatrix}
	-\Gamma_{21} C_qe+\Gamma_{22}\Delta\tilde p\\
	\Delta\tilde p
	\end{pmatrix}^\top \tilde M_2
	\begin{pmatrix}
	-\Gamma_{21} C_qe+\Gamma_{22}\Delta\tilde p\\
	\Delta\tilde p
	\end{pmatrix}\nonumber\\
	=&
	\begin{pmatrix}
	\Delta q\\
	\Delta p
	\end{pmatrix}^\top T_2^\top\tilde M_2 T_2
	\begin{pmatrix}
	\Delta q\\
	\Delta p
	\end{pmatrix}\geq 0,\nonumber
	\end{align}
	where  the last equality is from \eqref{eq:Mass2} in Assumption \ref{ass2}. Hence, $e^\top C_q^\top\Gamma_{21}^\top X_2^{-1}\Gamma_{21}C_qe-2e^\top C_q^\top\Gamma_{21}^\top X_2^{-1}\Gamma_{22}\Delta\tilde p+\Delta\tilde p^\top(\Gamma_{22}^\top X_2^{-1} \Gamma_{22}-Y_2^{-1})\Delta\tilde p\geq 0.$ 
	From \eqref{Mpf1}, the inequality above implies that $\kappa_1 \|e\|^2+\kappa_2\|e\|\|\Delta\tilde p\|-\kappa_3\|\Delta\tilde p\|^2\geq 0$
	where $\kappa_1=\lambda_{\max}(C_q^\top\Gamma_{21}^\top X_2^{-1}\Gamma_{21}C_q)$, $\kappa_2=2\|C_q^\top\Gamma_{21}^\top X_2^{-1}\Gamma_{22}\|$, $\kappa_3=\lambda_{\min}(Y_2^{-1}-\Gamma_{22}^\top X_2^{-1} \Gamma_{22})$. Clearly, $\kappa_1,\kappa_3>0$, $\kappa_2\geq 0$. Therefore, we have
	$\|\Delta\tilde p\|\leq \kappa \|e\|$ 
	where $\kappa:=(\kappa_2+\sqrt{\kappa_2^2+4\kappa_1\kappa_3})/2\kappa_3>0$.

	Since $\Delta k=\tilde{K}_1e+\tilde{K}_2\Delta \tilde p$ by \eqref{Delk}, we have 
	\begin{align}\label{delk}
	\|\Delta k\|  \leq \hat \kappa\|e\|
	\end{align}
	for all $x,e$, where $\hat{\kappa}=\|\tilde{K}_1\|+\|\tilde{K}_2\|\kappa>0$, which bounds $\|\Delta k\|/\|e\|$.

	4)
	Next, we analyse controller performance.
	From \eqref{L1} we have $R_3=\tilde K_1P_2$ and $R_4=\tilde K_2Y_2$ where $\tilde K_1,\tilde K_2$ are given in \eqref{tildeK}. Plugging $R_3,R_4$ into \eqref{Psi}, we have $\Psi_0=(\tilde A_2+B\tilde K_1)P_2+P_2(\tilde A_2+B\tilde K_1)^\top+\alpha_2P_2$, and the $(1,2)$ entry of $\Psi$ to be $(\tilde E_2+B\tilde K_2)Y_2$. Pre- and post-multiplying the inequality \eqref{LMI2} by the matrix $diag(I_n,Y_2^{-1},I_{n_w},I_{n_q})$, and then applying Schur's complement, we have 
	\begin{align}
	&\tilde\Psi+
	\begin{pmatrix}
	\psi_1 \\
	\varphi
	\end{pmatrix}^\top\tilde M_2
	\begin{pmatrix}
	\psi_1 \\
	\varphi
	\end{pmatrix}\preceq 0,\label{eqPsi1}
	\end{align}
	where  
	\begin{align*}
	\tilde\Psi&=
	\begin{pmatrix}
	\Psi_0& \tilde E_2+B\tilde K_2 & E_w \\
	*   & {\bf 0} &{\bf 0}  \\
	*  & *&-\mu_2I
	\end{pmatrix},
	\end{align*}
	and $\psi_1=(\Gamma_{21} C_q P_2 ,  \Gamma_{22} ,{\bf 0}_{n_q\times n_w})$, $\Psi_0$ is shown above, $\varphi$ is given in \eqref{varphi}.
	Let $P_3=P_2^{-1}$ and  pre- and post-multiply the inequality \eqref{eqPsi1} by $diag(P_3,I_{n_p},I_{n_w})$ and its transpose, respectively. This results in 
	\begin{align}
	&\hat\Psi+
	\begin{pmatrix}
	\psi_0 \\
	\varphi
	\end{pmatrix}^\top\tilde M_2
	\begin{pmatrix}
	\psi_0 \\
	\varphi
	\end{pmatrix}\preceq 0,\label{eqPsi2}
	\end{align}
	where $\psi_0=(\Gamma_{21} C_q,  \Gamma_{22} ,{\bf 0}_{n_q\times n_w})$ and
	\begin{align}
	\hat\Psi&=
	\begin{pmatrix}
	\hat\Psi_0& P_3(\tilde E_2+B\tilde K_2) & P_3E_w \\
	*  & {\bf 0} &{\bf 0}  \\
	* & *&-\mu_2I
	\end{pmatrix},\label{hatpsi}\\
	\hat\Psi_0&=P_3(\tilde A_2+B\tilde K_1)+(\tilde A_2+B\tilde K_1)^\top P_3+\alpha_2P_3.\label{hatpsi0}
	\end{align}
	
	Define $\xi_2=(x^\top,\tilde p_2^\top,w^\top)^\top$. Pre- and post-multiplying the inequality \eqref{eqPsi2} by $\xi_2^\top$ and $\xi_2$, respectively, we have 
	\begin{align}
	\xi_2^\top \hat\Psi\xi_2+\xi_2^\top
	\begin{pmatrix}
	\psi_0 \\
	\varphi
	\end{pmatrix}^\top\tilde M_2
	\begin{pmatrix}
	\psi_0 \\
	\varphi
	\end{pmatrix}\xi_2\leq 0.\label{ineqV21}
	\end{align}
	
	By \eqref{eq:Mass2} and \eqref{tilq}, we have
	$\begin{pmatrix}
	\tilde q_2\\
	\tilde p_2
	\end{pmatrix}^\top \tilde M_2 
	\begin{pmatrix}
	\tilde q_2\\
	\tilde p_2
	\end{pmatrix}\geq 0$. 
	Since $\tilde{q}_2=\Gamma_{21}q+\Gamma_{22}\tilde p_2=\Gamma_{21}C_qx+\Gamma_{22}\tilde p_2$,  $\begin{pmatrix}
	\tilde q_2\\
	\tilde p_2
	\end{pmatrix}=\begin{pmatrix}
	\psi_0 \\
	\varphi
	\end{pmatrix}\xi_2$, $\xi_2^\top
	\begin{pmatrix}
	\psi_0 \\
	\varphi
	\end{pmatrix}^\top\tilde M_2
	\begin{pmatrix}
	\psi_0 \\
	\varphi
	\end{pmatrix}\xi_2\geq 0.$ 
	Thus, $\xi_2^\top \hat\Psi\xi_2\leq 0$ from \eqref{ineqV21}, which is equivalent to $2x^\top P_3[(\tilde A_2+B\tilde K_1)x+(\tilde E_2+B\tilde K_2)\tilde p_2+E_ww]+\alpha_2 x^\top P_3x-\mu_2\|w\|^2\leq 0.$

	Let $V_2(x)=x^\top P_3x$. Then the derivative of $V_2$ along the trajectory of \eqref{dyn3-2} satisfies 
	\begin{align*}
	\dot V_2 &=2x^\top P_3[(\tilde A_2+B\tilde{K}_1)x+(\tilde E_2+B\tilde{K}_2)\tilde p_2-B\Delta k+E_ww]\\
	&\leq -\alpha_2 x^\top P_3x+\mu_2\|w\|^2+2\|P_3B\|\|x\|\|\Delta k\|.
	\end{align*}
	Recalling \eqref{delk}, we have 
	\begin{equation}
	\dot{V}_2 \leq -\alpha_2 x^\top P_3x+\mu_2\|w\|^2+\theta\|x\|\|e\|.\label{eqxthm1}
	\end{equation}
	where 
	$\theta=2\|P_3B\|\hat{\kappa}$.

	5) Finally, we prove that the  closed-loop system expressed by \eqref{dyn3-2} and \eqref{dyn4-2} is ISS with respect to $w$.
	Choose two constants $c_1,c_2$ as 
	$c_1=\alpha_1\lambda_m(P_1)/\lambda_M(P_1),\;c_2=\alpha_1\lambda_m(P_3)/\lambda_M(P_3).$ 
	Since $c_1>0,c_2>0$, we can choose two constants $\alpha_0>0,\beta_0>0$ such that
	$\alpha_0\!<\!\min\{c_1,c_2\}$, $\beta_0\!\geq\! \frac{\theta^2}{4\lambda_M(P_1)\lambda_M(P_3)(c_1-\alpha_0)(c_2-\alpha_0)}$. 
	Then, it is easy to check that the  matrix $P_0:= 
	\begin{pmatrix}
	\tilde P_0&\theta/2\\
	\theta/2&\hat P_0
	\end{pmatrix}$ 
	is negative semi-definite  where $\tilde P_0=-\alpha_2\lambda_m(P_3)+\alpha_0 \lambda_M(P_3)$ and $\hat P_0=\beta_0(-\alpha_1\lambda_m(P_1)+\alpha_0 \lambda_M(P_1))$. 
	Define a matrix $P$ as 
	$P=
	\begin{pmatrix}
	P_3&{\bf 0}\\
	{\bf 0}&\beta_0P_1
	\end{pmatrix}$. 
	Clearly, $P$ is positive definite. We can verify that the candidate Lyapunov function $V(x,e):=z^\top Pz$
	satisfies $V(x,e)=\beta_0V_1(e)+ V_2(x)$, and its derivative along the trajectory of \eqref{dyn3-2} and \eqref{dyn4-2} satisfies 
	\begin{align}
	\dot V+\alpha_0 V
	\leq& -\alpha_1 \beta_0 e^\top P_1e -\alpha_2 x^\top P_3x+\theta\|x\|\|e\|+\alpha_0 V\nonumber\\
	\leq&(\|x\|,\|e\|) P_0(\|x\|,\|e\|)^\top+(\mu_1\beta_0+\mu_2)\|w\|^2\nonumber\\
	\leq &(\mu_1\beta_0+\mu_2)\|w\|^2.\label{eqxethm1}
	\end{align}
	Therefore, the closed-loop system \eqref{dyn3-2} and \eqref{dyn4-2}, or equivalently \eqref{dyn3-1}, satisfies  \eqref{ISSineq} with $\mathcal{K}_\infty$ functions $\gamma(\|(x,e)\|)=-\alpha_0\lambda_{m}(P)\|(x,e)\|^2$ and $\chi(\|w\|)=(\mu_1\beta_0+\mu_2)\|w\|^2$.  This completes the proof.
\end{proof}

\begin{remark}
In the proof of Theorem \ref{thm1}, we first prove that the observer error $e$ is ISS w.r.t. $w$ (see  \eqref{eqethm1} in step 2), then the state $x$ is ISS w.r.t. $w$ and $e$ (see \eqref{eqxthm1} in step 4), and finally $(x,e)$ is ISS w.r.t. $w$ (see \eqref{eqxethm1} in step 5). 
This procedure of proving global stabilization is similar to the certainty equivalence proof  used in \cite{tsinias1993sontag}.
\end{remark}

\begin{remark}	
If $\alpha_1$ is fixed, then \eqref{LMI1} is an LMI in decision variables $\mu_1,P_1,R_1,R_2,X_1,Y_1$ that are used to determine observer gains $L_1,L_2$; if $\alpha_2$ is fixed, then \eqref{LMI2} is an LMI in decision variables $\mu_2,P_2,R_3,R_4,X_2,Y_2$ that are used to determine controller gains $K_1,K_2$. There is no coupling in decision variables in LMIs \eqref{LMI1}  and \eqref{LMI2}, implying a separation of the controller and observer designs.

\end{remark}

\begin{remark}\label{rem-recompute}
	The proof of Theorem \ref{thm1} indicates that
	larger $\alpha_1,\alpha_2$ result in a larger function $\gamma(\cdot)\in\mathcal{K}_\infty$ in \eqref{ISSineq}, which in turn indicates a faster convergence rate for the system \eqref{eqthm3-1}. 
	The convergence rate guarantee given in the proof of Theorem \ref{thm1}, $\alpha_0$, can be improved by finding a new ISS-Lyapunov function $V=z^\top Pz$ via  Proposition \ref{prp1}.
\end{remark}

\begin{remark}
	The condition $M_{24}\prec 0$ in Theorem \ref{thm1} is  used to prove that $\|\Delta k\|/\|e\|$ is bounded, which holds automatically when $p$ is globally Lipschitz. The condition $M_{24}\prec 0$ can be replaced by a growth condition on $p$ similar to Theorem 2 of \cite{arcak2001observer}. 
	Specifically, for a system described by \eqref{dyn1}-\eqref{eq:delQC} where $E_w=F_w={\bf 0}$, if all the conditions of Theorem \ref{thm1} but $M_{24}\prec0$ hold,  and there exist a $\mathcal{K}$ function $g_1$ and a non-decreasing function $g_2:[0,\infty)\rightarrow[0,\infty)$ such that $\|p(C_q(x+\Delta x))-p(C_qx)\|\leq g_1(\|\Delta x\|)\|C_qx\|$ 
	for all $x,\Delta x$ that satisfy $\|C_qx\|\geq g_2(\|\Delta x\|)$,  then 
	the feedback controller \eqref{input1} with $L_1,L_2,K_1,K_2$ given by \eqref{L1} renders the  closed-loop system \eqref{dyn3-1} globally  exponentially stable.

	The condition $M_{24}\prec0$ can be eliminated by using a  simpler form of $u$. Specifically, suppose that the observer-based controller $u$ has the form $u(t)=K_1\hat x(t)$, all the conditions of Theorem \ref{thm1} but $M_{24}\prec0$ hold with $R_4={\bf 0}$, and $L_1,L_2,K_1$ are given by \eqref{L1}, then the  closed-loop system \eqref{dyn3-1} is ISS w.r.t. $w$. In this case, the LMI \eqref{LMI2} is less likely to be satisfied by fixing $R_4={\bf 0}$.
	
\end{remark}

\subsection{Block Anti-Triangular Parameterization}
In this subsection, we consider a block anti-triangular parameterization of the $\delta$-MM for $p$. The following assumption on the parameterization of $M$ is given first.

	\begin{assumption}\label{ass3}
		There exist a set $\mathcal{N}$ of matrix pairs $(X,Y)$ with $X\in\R^{n_q\times n_p}$, $Y\in\R^{n_p\times n_p}$, and an invertible matrix $T\in\R^{(n_p+n_q)\times (n_p+n_q)}$, such that $M$ given below is a $\delta$-MM of $p$ for all $(X,Y)\in\mathcal{N}$:
		\begin{align}
		M=T^\top \tilde M T\;\mbox{where}\;\tilde M=
		\begin{pmatrix}
		{\bf 0}&X\\
		X^\top&Y
		\end{pmatrix}.\label{eq:Mass3}
		\end{align}  
	\end{assumption}

The following theorem  provides sufficient conditions for the design of matrices $L_1,L_2,K_1,K_2$ when the $\delta$-MM $M$ can be parameterized in a block anti-triangular manner. 

\begin{theorem}\label{thm2}
	Consider the system described by \eqref{dyn1}-\eqref{eq:delQC} with $p({\bf 0})={\bf 0}$.  Suppose that \\
		1) Assumption \ref{ass3} holds for some $T_1,\mathcal{N}_1$ and $T_2,\mathcal{N}_2$, respectively, where $T_1$ and $T_2$ are partitioned as in \eqref{eqT} and in \eqref{eqT2}, respectively, with $T_{14}$, $T_{24}$ invertible;\\
		2) 
		there exist positive constants $\alpha_1,\alpha_2,\mu_1,\mu_2$, matrices $R_1,R_2,R_3,R_4,X_1,Y_1.X_2,Y_2$, and symmetric and positive definite matrices $P_1,P_2$, such that $(X_1,Y_1)\in\mathcal{N}_1$, $(X_2,Y_2)\in\mathcal{N}_2$, and 
	\begin{align}
	&\Phi+\Upsilon_1^\top \tilde M_1\Upsilon_1+\Upsilon_2^\top \Upsilon_1+\Upsilon_1^\top \Upsilon_2\preceq 0,\label{LMI3}\\
	&\Psi+\Upsilon_3^\top \tilde M_2\Upsilon_3+\Upsilon_4^\top \Upsilon_3+\Upsilon_3^\top \Upsilon_4\preceq 0,\label{LMI4}\\
	&\Gamma_{12}^\top X_1+X_1^\top\Gamma_{12}+Y_1\prec 0,\label{LMI5}
	\end{align}
	where
	\begin{align}
	\tilde M_1&=
	\begin{pmatrix}
	{\bf 0}&X_1\\
	X_1^\top&Y_1
	\end{pmatrix},\;\tilde M_2=
	\begin{pmatrix}
	{\bf 0}&X_2\\
	X_2^\top&Y_2
	\end{pmatrix},\\
	\Psi&=
	\begin{pmatrix}
	\Psi_0& \tilde E_2+BR_4 & E_w \\
	*  & {\bf 0} &{\bf 0}  \\
	* & *&-\mu_2I
	\end{pmatrix},
	\end{align}
	\begin{align}\label{upsilon1}
	\begin{cases}
	\Upsilon_1=
	\begin{pmatrix}
	-\Gamma_{11} C_q& \Gamma_{12}& {\bf 0}_{n_q\times n_w} \\
	{\bf 0}_{n_p\times n_x}& I_{n_p}& {\bf 0}_{n_p\times n_w}
	\end{pmatrix},\\
	\Upsilon_2=
	\begin{pmatrix}
	{\bf 0}_{n_q\times n_x} & {\bf 0}_{n_q\times n_p} & {\bf 0}_{n_q\times n_w} \\
	-R_2 C& {\bf 0}_{n_p} & -R_2F_w \\
	\end{pmatrix},\\
	\Upsilon_3=
	\begin{pmatrix}
	{\bf 0}_{n_q\times n_x} & \Gamma_{22}& {\bf 0}_{n_q\times n_w} \\
	{\bf 0}_{n_p\times n_x} & I_{n_p}& {\bf 0}_{n_p\times n_w}  \\
	\end{pmatrix},\\
	\Upsilon_4=
	\begin{pmatrix}
	{\bf 0}_{n_q\times n_x} & {\bf 0}_{n_q\times n_p} & {\bf 0}_{n_q\times n_w} \\
	X_2\Gamma_{21}C_qP_2& {\bf 0}_{n_p} & {\bf 0}_{n_p\times n_w}  \\
	\end{pmatrix},
	\end{cases}
	\end{align}
	with  $\Gamma_{i1}(i=1,2)$ given in \eqref{gamma1}, $\Gamma_{i2}(i=1,2)$ given in \eqref{gamma2}, $\Phi$ given in \eqref{Phi},  $\Phi_0$  given in \eqref{Phi0},   $\Psi_0$ given in \eqref{Psi0}, $\varphi$ given in \eqref{varphi}, $\tilde A_i(i=1,2)$ given in \eqref{tilA} and $\tilde E_i(i=1,2)$ given in \eqref{tilE}. If $X_1$ has full row rank, then the  closed-loop system \eqref{dyn3-1} is ISS w.r.t. $w$ with
		$L_2$, $K_1$ given by \eqref{L1},  and $L_1,K_2$ given by
		\begin{align}\label{L1v2}
		L_1=\Gamma_{11}^{-1}X_1^{\dagger}R_2,\quad
		K_2=R_4 T_{24},
		\end{align} 
		where $X_1^{\dagger}$ is the right inverse of $X_1$.
\end{theorem}
\begin{proof}
	As shown in \eqref{dyn3-2} and \eqref{dyn4-2}, dynamics of the closed-loop system under transformations can be described as
	\begin{align*}
	\dot x&=(\tilde A_2+B\tilde{K}_1)x+(\tilde E_2+B\tilde{K}_2)\tilde p_2-B\Delta k+E_ww,\\
	\dot e&=(\tilde A_1+\tilde L_2C)e-\tilde E_1\delta\tilde p_1+(E_w+\tilde L_2F_w)w,
	\end{align*}
	where $\tilde{K}_1,\tilde{K}_2$ are given in \eqref{tildeK}, $\Delta k$ is given in \eqref{Delk},  $\tilde p_2$ is given in \eqref{tilq}, and $\tilde L_2$ is given in \eqref{tilL2}.
		From \eqref{L1} and \eqref{L1v2}, we have $R_1=P_1\tilde L_2$ and $R_2=X_1\Gamma_{11}L_1$. 
	We claim that \eqref{LMI3} is equivalent to 
	\begin{align}
	&\Phi+Q_1^\top\tilde M_1Q_1\preceq 0\label{ineqV3}
	\end{align}
	and  \eqref{LMI4} is equivalent to 
	\begin{align}
	&\Psi+Q_2^\top\tilde M_2Q_2\preceq 0\label{ineqV5}
	\end{align}
	where
	\begin{align*}
	Q_1&=
	\begin{pmatrix}
	-\Gamma_{11} (C_q+L_1C)& \Gamma_{12}& -\Gamma_{11}L_1F_w\\
	{\bf 0}_{n_p\times n_x}& I_{n_p}& {\bf 0}_{n_p\times n_w} 
	\end{pmatrix},\\
	Q_2&=
	\begin{pmatrix}
	\Gamma_{21} C_qP_2&  \Gamma_{22} &{\bf 0}_{n_q\times n_w} \\
	{\bf 0}_{n_p\times n_x}& I_{n_p}& {\bf 0}_{n_p\times n_w} 
	\end{pmatrix}.
	\end{align*}
	
	Indeed, $Q_1$ can be written as $Q_1=\Upsilon_1+\hat \Upsilon_2$ where $\Upsilon_1$ is given in \eqref{upsilon1} and 
	\begin{align*}
	\hat\Upsilon_2&=
	\begin{pmatrix}
	-\Gamma_{11}L_1 C& {\bf 0}_{n_q\times n_p}& -\Gamma_{11}L_1F_w \\
	{\bf 0}_{n_p\times n_x}& {\bf 0}_{n_p\times n_p}& {\bf 0}_{n_p\times n_w}
	\end{pmatrix}.
	\end{align*}
	It is easy to verify that $\Upsilon_2=\tilde M_1 \hat\Upsilon_2$ and $\hat\Upsilon_2^\top \tilde M_1\hat\Upsilon_2={\bf 0}$. 
	Therefore,
	\begin{align*}
	Q_1^\top\tilde M_1Q_1&=(\Upsilon_1+\hat\Upsilon_2)^\top\tilde M_1(\Upsilon_1+\hat\Upsilon_2)\\
	&=\Upsilon_1^\top \tilde M_1\Upsilon_1+\hat\Upsilon_2^\top \tilde M_1\Upsilon_1+\Upsilon_1^\top \tilde M_1\hat\Upsilon_2+\hat\Upsilon_2^\top \tilde M_1\hat\Upsilon_2\\
	&=\Upsilon_1^\top \tilde M_1\Upsilon_1+\Upsilon_2^\top \Upsilon_1+\Upsilon_1^\top \Upsilon_2.
	\end{align*}
	Similarly, $Q_2$ can be written as $Q_2=\Upsilon_3+\hat \Upsilon_4$ where $\Upsilon_3$ is given in \eqref{upsilon1} and 
	\begin{align*}
	\hat\Upsilon_4&=
	\begin{pmatrix}
	\Gamma_{21}C_qP_2& {\bf 0}_{n_q\times n_p}& {\bf 0}_{n_q\times n_w} \\
	{\bf 0}_{n_p\times n_x}& {\bf 0}_{n_p\times n_p}& {\bf 0}_{n_p\times n_w}
	\end{pmatrix}.
	\end{align*}
	It is easy to verify that $\Upsilon_4=\tilde M_2 \hat\Upsilon_4$ and $\hat\Upsilon_4^\top \tilde M_2\hat\Upsilon_4={\bf 0}$.
	Therefore, $Q_2^\top\tilde M_2Q_2=(\Upsilon_3+\hat\Upsilon_4)^\top\tilde M_2(\Upsilon_3+\hat\Upsilon_4)
	=\Upsilon_3^\top \tilde M_2\Upsilon_3+\Upsilon_3^\top \Upsilon_4+\Upsilon_4^\top \Upsilon_3.$ 
	Hence, our claim is proved.

	Plugging $R_1$ into $\Phi_0$ and $\Phi$, we have $\Phi_0=P_1(\tilde A_1+\tilde L_2C)+(\tilde A_1+\tilde L_2C)^\top P_1+\alpha_1 P_1$, and the $(1,3)$ entry of $\Phi$ is $P_1(E_w+\tilde L_2F_w)$.
	Define $\xi_1=(e^\top,\delta\tilde p_1^\top,w^\top)^\top$. Pre- and post-multiplying \eqref{ineqV3} by $\xi_1^\top$ and $\xi_1$, respectively, we have $\xi_1^\top \Phi\xi_1+\xi_1^\top Q_1^\top\tilde M_1Q_1\xi_1\leq 0.$
	Since $Q_1\xi_1=
	\begin{pmatrix}
	\delta\tilde q_1\\
	\delta\tilde p_1
	\end{pmatrix}=T_1
	\begin{pmatrix}
	\delta q \\
	\delta p
	\end{pmatrix}$ and $M_1$ satisfies Assumption \ref{ass3}, 
	we have $\xi_1^\top Q_1^\top\tilde M_1Q_1\xi_1\geq 0$, 
	which implies that $\xi_1^\top \Phi\xi_1\leq 0$. Hence, $2e^\top P_1[(\tilde A+\tilde L_2\tilde C)e-\tilde E_1 \delta\tilde p_1+(E_w+\tilde L_2F_w)w]+\alpha_1 e^\top P_1e-\mu_1\|w\|^2\leq 0.$ 
	Define $V_1(e)=e^\top P_1e$. Then, 
	we have 
	\begin{align*}
	\dot V_1\leq -\alpha_1 e^\top P_1e+\mu_1\|w\|^2.
	\end{align*}
	
	Define $\Delta q=C_q\hat x-C_qx$ and $\Delta \tilde q:=\tilde q_1(C_q\hat x)-\tilde q_1(C_qx)$.  Then, $\Delta q=-C_qe$ and $\Delta \tilde q=-\Gamma_{11}C_qe+\Gamma_{12}\Delta\tilde p$ where $\Delta \tilde{p}=\tilde{p}_1 (C_q\hat{x})-\tilde{p}_1(C_qx)$. Define  $\zeta=(e^\top,\Delta\tilde p^\top)^\top$. Therefore,
	\begin{align*}
	&\zeta^\top
	\begin{pmatrix}
	-\Gamma_{11}C_q & \Gamma_{12} \\
	{\bf 0} & I
	\end{pmatrix}^\top\tilde M_1
	\begin{pmatrix}
	-\Gamma_{11}C_q & \Gamma_{12} \\
	{\bf 0} & I
	\end{pmatrix}\zeta\nonumber\\
	=&
	\begin{pmatrix}
	\Delta q\\
	\Delta p
	\end{pmatrix}^\top T_1^\top\tilde M_1T_1
	\begin{pmatrix}
	\Delta q\\
	\Delta p
	\end{pmatrix}\geq 0,
	\end{align*}
	where the last equality is from Assumption \ref{ass3}. Hence, $-2e^\top C_q^\top\Gamma_{11}^\top X_1 \Delta\tilde p+\Delta\tilde p^\top(\Gamma_{12}^\top X_1+X_1^\top \Gamma_{12}+Y_1)\Delta\tilde p\geq 0.$ 
	From \eqref{LMI5}, the inequality above implies that
	$\kappa_1 \|e\|\|\Delta\tilde p\|-\kappa_2\|\Delta\tilde p\|^2\geq 0$,
	where $\kappa_1=2\|C_q^\top\Gamma_{11}^\top X_1\|$ and $\kappa_2=-\lambda_{\min}(\Gamma_{12}^\top X_1+X_1^\top \Gamma_{12}+Y_1)$. Noticing that $\kappa_1\geq 0,\kappa_2>0$, we have $\|\Delta\tilde p\|\leq \frac{\kappa_1}{\kappa_2} \|e\|$. 
	Noting that $\Delta k=\hat{K}_1e+\hat{K}_2\Delta \tilde p$ with $\Delta \tilde{p}=\tilde{p}_1 (C_q\hat{x})-\tilde{p}_1(C_qx)$, $\hat{K}_1= K_1 -K_2T_{14}^{-1}T_{13}C_q$, $\hat{K}_2 = K_2T_{14}^{-1}$, we have 
	$\|\Delta k\|  \leq \hat \kappa\|e\|$ for all $x,e$, where $\hat{\kappa}=\|\tilde{K}_1\|+\|\tilde{K}_2\|\kappa_1/\kappa_2\geq 0$.

	From \eqref{L1} and \eqref{L1v2}  we have $R_3=\tilde K_1P_2$ and $R_4=\tilde K_2$ where $\tilde K_1,\tilde K_2$ are defined in \eqref{tildeK}. Plugging $R_3,R_4$ into $\Psi_0$ and $\Psi$, we have $\Psi_0=(\tilde A_2+B\tilde K_1)P_2+P_2(\tilde A_2+B\tilde K_1)^\top+\alpha_2P_2$, and the $(1,2)$ entry of $\Psi$ is $\tilde E_2+B\tilde K_2$. 
	Let $P_3=P_2^{-1}$ and  pre- and post-multiply \eqref{ineqV5} by $diag(P_3,I_{n_p},I_{n_w})$ and its transpose, respectively. This results in 
	\begin{align}
	&\hat\Psi+Q_3^\top\tilde M_2Q_3\preceq 0,\label{eqPsi3}
	\end{align}
	where $Q_3=
	\begin{pmatrix}\Gamma_{21} C_q&  \Gamma_{22} &{\bf 0}_{n_q\times n_w}\\
	{\bf 0}_{n_p\times n_x}& I_{n_p}& {\bf 0}_{n_p\times n_w} \end{pmatrix}$
	and $\hat\Psi$ is given in \eqref{hatpsi} with $\hat\Psi_0$ given in \eqref{hatpsi0}. Define $\xi_2=(x^\top,\tilde p_2^\top,w^\top)^\top$. Pre- and post-multiplying \eqref{eqPsi3} by $\xi_2^\top$ and $\xi_2$, respectively, we have $\xi_2^\top \hat\Psi\xi_2+\xi_2^\top Q_3^\top\tilde M_2Q_3\xi_2\leq 0$. 
	Since  $Q_3\xi_2=
	\begin{pmatrix}
	\tilde q_2\\
	\tilde p_2
	\end{pmatrix}=T_2
	\begin{pmatrix}
	q \\
	p
	\end{pmatrix}$ and $M_2$ satisfies Assumption \ref{ass3}, it follows that $\xi_2^\top Q_3^\top\tilde M_2Q_3\xi_2\geq 0$. 
	Hence, we have $\xi_2^\top \hat\Psi\xi_2\leq 0$, which is equivalent to $2x^\top P_3[(\tilde A_2+B\tilde K_1)x+(\tilde E_2+B\tilde K_2)\tilde p_2+E_ww]+\alpha_2 x^\top P_3x-\mu_2\|w\|^2\leq 0.$ 
	Let $V_2(x)=x^\top P_3x$. Then, we have $\dot V_2=2x^\top P_3[(\tilde A_2+B\tilde{K}_1)x+(\tilde E_2+B\tilde{K}_2)\tilde p_2-B\Delta k+E_ww]\leq -\alpha_2 x^\top P_3x+\mu_2\|w\|^2+2\|P_3B\|\|x\|\|\Delta k\|.$ 
	Recalling $\|\Delta k\|  \leq \hat \kappa\|e\|$, we have
	\begin{align*}
	\dot{V}_2 \leq -\alpha_2 x^\top P_3x+\mu_2\|w\|^2+\theta\|x\|\|e\|
	\end{align*}
	where $\theta=2\|P_3B\|\hat{\kappa}$. The rest of the proof proceeds as that given in part 5) of the proof of Theorem \ref{thm1}.
\end{proof}

\begin{remark}
Inequality \eqref{LMI3} is an LMI in decision variables $\mu_1,P_1,R_1,R_2$ when $\alpha_1$ is fixed, \eqref{LMI4} is an LMI in  decision variables $\mu_2,P_2,R_3,R_4$ when  $\alpha_2$ and $X_2$ are fixed, and \eqref{LMI5} is an LMI in  decision variables $X_1,Y_1$. Hence, we can fix $\alpha_1,\alpha_2,X_2$ and solve for  \eqref{LMI3}-\eqref{LMI5}.
When $L_1,L_2,K_1,K_2$ are obtained, a re-computation for $P,\alpha_0,\mu$ using Proposition \ref{prp1} may result in a better convergence rate guarantee. 
\end{remark}

\begin{remark}The LMIs \eqref{LMI1} and \eqref{LMI2}  both have dimensions $(n_x+n_p+n_q+n_w)\times (n_x+n_p+n_q+n_w)$, the LMIs \eqref{LMI3} and \eqref{LMI4}  both have dimensions $(n_x+n_p+n_w)\times (n_x+n_p+n_w)$, and the LMI \eqref{LMI5} has dimension $n_p\times n_p$. These LMIs can be solved reliably and efficiently by the interior point method algorithms of convex optimization with a polynomial-time complexity. Exploring for what class of systems these LMIs are guaranteed to be feasible (i.e., analytical verification of feasibility) is still under our investigation. Furthermore, these LMI conditions might be conservative compared with specific results that focus on certain special nonlinearities such as the globally  Lipschitz nonlinearity.

\end{remark}

\section{Event-triggered Control Design}\label{sec:trigger}

In this section, we  discuss event-triggering mechanisms (ETMs) within the observer-based controller designed in the preceding section for the system described by \eqref{dyn1}-\eqref{eq:delQC} where the nonlinearity $p$ is assumed to be globally Lipschitz. For certain  incrementally quadratic nonlinearities that imply the global Lipschitzness (such as the incremental sector bounded nonlinearity and the nonlinearities with Jacobians in polytopes \cite{accikmecse2011observers}), using their corresponding incremental matrix characterizations, instead of the matrix characterizations for global Lipschitzness, makes the associated LMIs in the design procedure less conservative, while benefiting from having the Lipschitz property needed for the upcoming ETM-related results to hold.

\subsection{Configuration I: The  Controller  Channel Is Implemented By ETM}

\begin{figure}[!hbt]
	\centering
	\includegraphics[height=3.5cm]{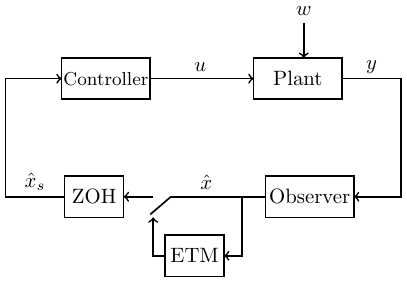}
	\caption{Configuration where the ETM is implemented in the controller channel.}\label{figConfig1}
\end{figure}

In this subsection, we discuss the configuration shown in Figure \ref{figConfig1} where  the plant is described by \eqref{dyn1}-\eqref{eq:delQC},  the observer is given in \eqref{obser1}, the continuous-time feedback controller is given in \eqref{input1}, and the ETM only has the information of $\hat x$, the state of the observer. We will assume that $\|w\|_\infty\leq \omega_0$ where $\omega_0$ is a positive constant indicating the bound of the disturbance in this subsection and the next subsection.

The feedback controller $u(t)$ is implemented by an ETM such that it is only updated at certain triggering time instances $t_1,t_2,...$ where $t_k<t_{k+1}$ for any $k\geq 0$ and kept constant during consecutive time instances. Define $t_0=0$ and the piecewise constant signal $\hat x_s$ as
\begin{align*}
\hat x_s(t)=\hat x(t_k),\;\forall t\in[t_k,t_{k+1}).
\end{align*} 
Then the control input $u(t)$ is given by 
\begin{align}\label{inputform2}
u(t)=K_1\hat x_s(t)+K_2 p(C_q\hat x_s(t))
\end{align}  
where $K_1,K_2$ are matrices to be designed.
The input $u(t)$ has the same form as that in \eqref{input1}, but it is updated  at triggering time instances  $t_1,t_2,\dots$, which are determined by the following type of triggering rule:
\begin{align}
&t_{k+1}=\inf\{t\mid t\geq t_k+\tau,\;\|\hat x_e(t)\|> \sigma \|\hat x(t)\|+\epsilon\}\label{triggercon1}
\end{align}
where $\hat x_e$ is defined as $\hat x_e(t)=\hat x_s(t)-\hat x(t)$ and $\tau,\sigma,\epsilon$ are all positive numbers to be specified. The time-updating rule 
\eqref{triggercon1} guarantees that the inter-execution times $\{t_{k+1}-t_{k}\}$ are lower bounded by the built-in positive constant $\tau$, which means that Zeno phenomenon (i.e., infinite executions happen in a finite amount of time) will not occur \cite{tabuada2007event}.

\begin{remark}
	The triggering rule  \eqref{triggercon1} only depends on the information of $\hat x$ and $\hat x_e$, which are available from the  proposed observer. This triggering rule  
	is a combination of a mixed ETM and a time regularization technique. There are several motivations for choosing this type of rule. It was known that even in the absence of disturbances, inter-execution times of many ETMs converge to zero for output-based control configurations \cite{donkers2012output}.  To exclude the Zeno phenomenon, time regularization or periodic event-triggered control, which enforces a built-in lower bound for  inter-execution times, has been utilized in recent works on observer-based ETMs \cite{abdelrahim2017robust,tallapragada2012event,dolk2017output}. Furthermore,  mixed ETM is known to be robust to external disturbances or measurement noise, while  relative ETM and  absolute ETM have zero robustness to  disturbance/noise \cite{borgers2014event}. Additionally,  an event-triggering rule with  time regularization can benefit from using mixed ETMs in terms of the number of events that are generated (e.g., see Example 3 in \cite{borgers2014event}). 
\end{remark}

The closed-loop system that combines system \eqref{dyn1}-\eqref{eq:delQC}, observer \eqref{obser1} and  event-triggered controller  \eqref{inputform2} is expressed compactly as
\begin{align}\label{closedtrigger}
\dot z=A_cz+H_1p+H_2\delta p+H_3\delta \hat p+H_4w+H_5\hat x_e
\end{align}
where $\delta p,\delta q$ are given in \eqref{delp}, $A_c$ is given in \eqref{eqA}, $H_1,H_2,H_3,H_4$ are given in \eqref{H1}, and
\begin{align}
\delta \hat p&=p(C_q\hat x_s)-p(C_qx),\label{deltahatp}\\
H_5&=
\begin{pmatrix}
BK_1\\
{\bf 0}
\end{pmatrix}.\label{H5}
\end{align}

\begin{theorem}\label{thm3}
		Consider the configuration shown in Figure \ref{figConfig1}  where the plant is described by \eqref{dyn1}-\eqref{eq:delQC} with $p({\bf 0})={\bf 0}$ and $\|w\|_\infty\leq \omega_0$ with $\omega_0$ a positive number. Suppose that  there exists $\ell>0$ such that $\|p(r)-p(s)\|\leq \ell\|r-s\|$ for any $r,s$. 
		Suppose that 
		there exist positive numbers $\alpha_0>0,\mu>0$, and matrices $P\succ 0,K_1,K_2,L_1,L_2$ such that the closed-loop system \eqref{eqthm3-1}  with  controller \eqref{input1} and  observer  \eqref{obser1}  satisfies  $\dot V\leq -\alpha_0 V+ \mu\|w\|^2$ where $V=z^\top Pz$. 
		Choose any $\epsilon>0$ and 
		\begin{align}
		\sigma&=\frac{\varrho\alpha_0\lambda_m(P)}{2\sqrt{2}s}>0\label{eqsigma}
		\end{align}
		where $0<\varrho<1$ and $s=\|PH_5\|+\ell\|PH_3\|\|C_q\|$. 
		Choose $\tau>0$ as the solution to  the equation $\phi(\tau)=1$ where $\phi$ is the solution of the following ODE: 
		\begin{align*}
		&\dot \phi=\sqrt{2}(\eta_4+\eta_2\phi)(1+\sigma\phi),\;\phi(0) = 0,
		\end{align*}
		with 
		\begin{align}\label{eta1}
		\begin{cases}
		\eta_1=\|A_c\|+\ell\sqrt{b_1^2+b_2^2},\\
		\eta_2=\|H_5\|+\ell\|H_3\|\|C_q\|,\\
		\eta_3=\ell\|H_2\|\|L_1F_w\|+\|H_4\|,\\
		\eta_4=\frac{\eta_1}{\sqrt{2}\sigma}+\frac{\eta_3\omega_0}{\epsilon},\\
		b_1=\|H_1\|\|C_q\|,\\
		b_2=\|H_2\|\|C_q+L_1C\|+\|H_3\|\|C_q\|.
		\end{cases}
		\end{align}
		Then,  
		the closed-loop system \eqref{closedtrigger} that implements the triggering rule \eqref{triggercon1} is ISpS w.r.t. $w$. 
\end{theorem}
\begin{proof}
	Since the derivative of $V$ along the trajectory of the closed-loop system \eqref{eqthm3-1} satisfies $\dot V\leq -\alpha_0 V+ \mu\|w\|^2$, the derivative of $V$ along the trajectory of the closed-loop system \eqref{closedtrigger} satisfies $\dot V\leq -\alpha_0 V+ \mu\|w\|^2+2z^\top P[H_5\hat x_e+H_3(\delta \hat p-\Delta p)]
	\leq -\alpha_0\lambda_{m}(P)\|z\|^2+ \mu\|w\|^2+2\|z\|(\|PH_5\|\|\hat x_e\|
	+\|PH_3\|\|\delta \hat p-\Delta p\|).$ 
	Clearly,  $\|\delta \hat p-\Delta p\|=\|p(C_q\hat x_s)-p(C_q\hat x)\|\leq \ell\|C_q(\hat x_s-\hat x)\|\leq \ell\|C_q\|\|\hat x_e\|$. Then, we have
	\begin{align}
	\dot V&\leq -\alpha_0\lambda_{m}(P)\|z\|^2+ \mu\|w\|^2+2s\|z\|\|\hat x_e\|\nonumber\\
	&\leq -(1-\varrho)\alpha_0\lambda_{m}(P)\|z\|^2 + \mu\|w\|^2\nonumber\\
	&\quad \quad + \|z\|\Big[2s\|\hat x_e\|-\varrho\alpha_0\lambda_{m}(P)\|z\|\Big].\label{pf3dotv}
	\end{align} 
	For any $x,e$, we have $\|z\|=\sqrt{\|x\|^2+\|x-\hat x\|^2}=\sqrt{\|\hat x\|^2+2\|x\|^2-2x^\top\hat x}\geq \|\hat x\|/\sqrt{2}$, meaning that $\|\hat x\|\leq \sqrt{2}\|z\|$. Therefore, the condition 
	\begin{align}\label{triggercon}
	\|\hat x_e\|\leq \sigma \|\hat x\|+\epsilon
	\end{align}
	implies 
	\begin{align}
	\|\hat x_e\|\leq \sqrt{2}\sigma \|z\|+\epsilon,\label{pf3dotvcon}
	\end{align}
	which is equivalent to the inequality 
	$2s\|\hat x_e\|-\varrho\alpha_0\lambda_{m}(P)\|z\|\leq 2s\epsilon.$

	Choose a constant $c$ such that $0<c<(1-\varrho)\alpha_0\lambda_{m}(P)$. Then, as long as  \eqref{pf3dotvcon} holds, from \eqref{pf3dotv} we have
	\begin{align}
	\dot V
	&\leq-[(1-\varrho)\alpha_0\lambda_{m}(P)-c]\|z\|^2+\mu\|w\|^2+ \frac{s^2\epsilon^2}{c}.\label{triggerLayu}
	\end{align}

	Recalling that $p({\bf 0})={\bf 0}$ and $\ell$ is the Lipschitz constant of $p$, we have $\|p\|\leq \ell\|C_q\|\|x\|$, $\|\delta p\|\leq \ell \|\delta q\|\leq \ell(\|C_q+L_1C\|\|e\|+\|L_1F_w\|\|w\|)$, and $\|\delta \hat p\|\leq \ell \|C_q\|(\|\hat x_e\|+\|e\|)$. 
	Therefore, from \eqref{closedtrigger} we have 
	\begin{align*}
	\|\dot z\|
	&\leq \|A_c\|\|z\|+\|H_5\|\|\hat x_e\|+\ell\|H_1\|\|C_q\|\|x\|\\
	&\quad\quad +\ell\|H_2\|(\|C_q+L_1C\|\|e\|+\|L_1F_w\|\|w\|)\\
	&\quad\quad  +\ell\|H_3\|\|C_q\|(\|\hat x_e\|+\|e\|)+\|H_4\|\|w\|\\
	&\leq \eta_1\|z\|+\eta_2\|\hat x_e\|+\eta_3\|w\|
	\end{align*}
	where the second inequality follows from Cauchy's inequality  $b_1\|x\|+b_2\|e\|\leq \sqrt{b_1^2+b_2^2}\|z\|$. 
	
	Because $\|\dot z\|=\sqrt{\|\dot x\|^2+\|\dot e\|^2}=\sqrt{\|\dot x\|^2+\|\dot x-\dot{\hat{x}}\|^2}=\sqrt{\|\dot {\hat{x}}\|^2+2\|\dot x\|^2-2\dot x^\top\dot{\hat{x}}}\geq \|\dot{\hat{x}}\|/\sqrt{2}$  
	and $\|\dot{\hat{x}}_e\|=\|\dot{\hat{x}}\|$, we have $\|\dot{\hat{x}}_e\|\leq \sqrt{2}\|\dot z\|$.

	Let $v(t) = \frac{\|\hat x_e(t)\|}{\sqrt{2}\sigma\|z(t)\|+\epsilon}$. 
	Then for any $h>0$,
	\begin{align}
	&v(t+h)-v(t)=
	\frac{\|\hat x_e(t+h)\|}{\sqrt{2}\sigma\|z(t+h)\|+\epsilon}-\frac{\|\hat x_e(t)\|}{\sqrt{2}\sigma\|z(t)\|+\epsilon}
	\nonumber\\
	&=
	\frac{ \|\hat x_e(t\!+\!h)\|(\sqrt{2}\sigma\|z(t)\|\!+\!\epsilon)\!-\!\|\hat x_e(t)\|(\sqrt{2}\sigma\|z(t\!+\!h)\|\!+\!\epsilon)}{(\sqrt{2}\sigma\|z(t\!+\!h)\|\!+\!\epsilon)(\sqrt{2}\sigma\|z(t)\|\!+\!\epsilon)}
	\nonumber\\
	&=
	\frac{( \|\hat x_e(t+h) -\|\hat x_e(t)\|)(\sqrt{2}\sigma\|z(t)\|+\epsilon)}{(\sqrt{2}\sigma\|z(t+h)\|+\epsilon)(\sqrt{2}\sigma\|z(t)\|+\epsilon)}
	\nonumber\\
	&\quad -\frac{\sqrt{2}\sigma \|\hat x_e(t)\|(\|z(t+h)\|-\|z(t)\|)}{(\sqrt{2}\sigma\|z(t+h)\|+\epsilon)(\sqrt{2}\sigma\|z(t)\|+\epsilon)}
	\nonumber
	\end{align}
	and hence
	\begin{align}
	D^+v(t) &= \lim\sup_{h \rightarrow 0^+}\frac{v(t+h) -v(t)}{h}
	\nonumber\\
	&=\frac{D^+\|\hat{x}_e(t)\|}{\sqrt{2}\sigma\|z(t)\|+\epsilon}
	-\frac{\sqrt{2}\sigma\|\hat{x}_e(t)\|D^+\|z(t)\|}{(\sqrt{2}\sigma\|z(t)\|+\epsilon)^2}.
	\label{eq:D+v}
	\end{align}
	When $z(t) \neq0$, $D^+\|z(t)\| =\frac{z(t)^T\dot{z}(t)}{\|z(t)\|}$ and therefore $\lvert D^+\|z(t)\|\rvert\le \|\dot{z}(t)\|$.
	When $z(t)= 0$, $D^+\|z(t)\| = \lim\sup_{h \rightarrow 0^+} \frac{ \|z(t+h)\| -\|z(t)\| }{h}= \lim\sup_{h \rightarrow 0^+} \|\frac{z(t+h)}{h} \|=\| \dot{z}(t)\|$. 
	Thus, in all cases $\lvert D^+\|z(t)\|\rvert\le \|\dot{z}(t)\|$.
	Similarly, $\lvert D^+\|\hat{x}_e(t)\|\rvert\le \|\dot{\hat x}_e(t)\|$.
	Dropping the argument $t$,
	it now follows from \eqref{eq:D+v} that
	\begin{align*}
	D^+v &\le 
	\frac{\|\dot{\hat x}_e\|}{\sqrt{2}\sigma\|z\|+\epsilon}
	+\frac{\sqrt{2}\sigma\|\hat{x}_e\|\|\dot{z}\|}{(\sqrt{2}\sigma\|z\|+\epsilon)^2}
	\nonumber\\
	&\leq  \frac{\sqrt{2}\|\dot{z}\|}{\sqrt{2}\sigma\|z\|+\epsilon}+\frac{\sqrt{2}\sigma\|\hat x_e\|\|\dot z\|}{(\sqrt{2}\sigma\|z\|+\epsilon)^2}\nonumber\\
	&= \frac{\sqrt{2}\|\dot z\|}{\sqrt{2}\sigma\|z\|+\epsilon}(1+\frac{\sigma\|\hat x_e\|}{\sqrt{2}\sigma\|z\|+\epsilon})
	\nonumber\\
	&\leq  \sqrt{2}(\eta_4+\eta_2\frac{\|\hat x_e\|}{\sqrt{2}\sigma\|z\|+\epsilon})
	(1+\sigma\frac{\|\hat x_e\|}{\sqrt{2}\sigma\|z\|+\epsilon})
	\nonumber\\
	&=\sqrt{2}\left(\eta_4 +\eta_2v\right)\left(1 + \sigma v\right)
	\end{align*}
	where the following facts are used to derive the last inequality: 
	\begin{align*}
	\frac{\eta_1\|z\|}{\sqrt{2}\sigma\|z\|+ \epsilon}\leq \frac{\eta_1}{\sqrt{2}\sigma},\;\frac{\eta_3\|w\|}{\sqrt{2}\sigma\|z\|+\epsilon}\leq \frac{\eta_3\omega_0}{\epsilon}.
	\end{align*}
	
	Since $v(t_k) =0$, it now follows from  the comparison lemma  that $v(t)\leq \phi(t-t_k)$. 
	Since the time it takes for $v$ to evolve from $0$ to $1$ is lower bounded by $\tau$,  \eqref{pf3dotvcon} holds during the time interval $[t_k,t_k+\tau]$. 
	For any $k\geq 0$, if $t_{k+1}=t_k+\tau$, then \eqref{pf3dotvcon} holds during the interval  $[t_k.t_{k+1})$ as shown above; if $t_{k+1}>t_k+\tau$, then, during the interval $[t_k+\tau, t_{k+1})$, condition \eqref{triggercon} holds, which implies that  \eqref{pf3dotvcon} holds. Therefore, \eqref{pf3dotvcon}  holds during any interval $[t_k, t_{k+1})$ for any $k\geq 0$, i.e., it holds for any $t\geq 0$. 
	Since satisfaction of \eqref{pf3dotvcon} implies the inequality \eqref{triggerLayu}, we conclude that the function $V$ is an ISpS-Lyapunov function since it satisfies \eqref{ISSineq} for any $t\geq 0$ with $\gamma(\|z\|)=[(1-\varrho)\alpha_0\lambda_{m}(P)-c]\|z\|^2\in\mathcal{K}_\infty$, $\chi(\|w\|)=\mu\|w\|^2\in\mathcal{K}$ and $d=s^2\epsilon^2/c>0$. The conclusion follows by Proposition \ref{proISS}.
\end{proof}

	\begin{remark}\label{remarktradeoff}
		In the proof  of Theorem \ref{thm3}, the equation for $\tau$ is given explicitly, and any $\tau'\in(0,\tau]$ also makes the proof valid. 
		The parameter $\epsilon$ can be chosen arbitrarily, but there are trade-offs in choosing $\epsilon$: on one hand, the value of $d$ in the inequality \eqref{ISSineq} or \eqref{ISSineq2}   increases as $\epsilon$ increases, meaning that the ultimate bound for $x$ increases as $\epsilon$ increases; on the other hand, the explicit equation of $\tau$ depends on $\epsilon$, with $\tau$ decreasing to $0$ when $\epsilon$ approaches  $0$.  
		Hence, parameters in the triggering rule should be  chosen appropriately to balance the execution times and the performance.  Finding the maximal
		lower-bound of the inter-execution times is an interesting and challenging  problem that will
		be investigated in our future work.
	\end{remark}

\subsection{Configuration II:  The Controller  and   Observer Channels  Are Both Implemented By ETMs}
\begin{figure}[!hbt]
	\centering
	\includegraphics[height=3.5cm]{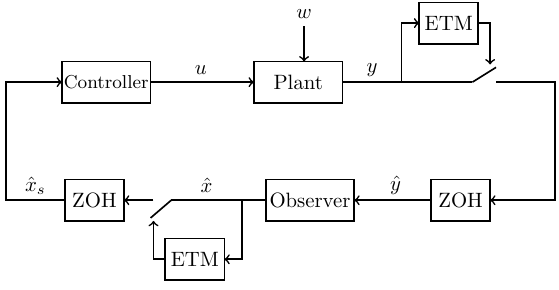}
	\caption{Configuration where ETMs are implemented in both the observer and  controller channels asynchronously}\label{figConfig2}
\end{figure}

In this subsection, we discuss the configuration shown in Figure \ref{figConfig2} where the ETM for the output is triggered by the information of $y$ and the ETM for the input is triggered by the information of $\hat x$, in an asynchronous manner.

Consider a system described by \eqref{dyn1}-\eqref{eq:delQC}.
The observer in the configuration of Figure \ref{figConfig2}  only has  sampled information $y_s(t)$ of the output $y(t)$ where $y_s(t)$ is updated at time instances $t_1^y,t_2^y,...$ by
\begin{align*}
y_s(t)=y(t_k^y),\;\forall t\in[t_k^y,t^y_{k+1}).
\end{align*} 
Here, $t_0^y=0$ and the triggering times $t_1^y,t_2^y,\dots$ are determined by the following triggering rule:
\begin{align}
&t_{k+1}^y=\inf\{t\mid t\geq t_k^y+\tau_y,\;\|y_e(t)\|> \sigma_y \|y(t)\|+\epsilon_y\}\label{triggeroutput1}
\end{align}
where $y_e(t)=y_s(t)-y(t)$ and $\tau_y,\sigma_y,\epsilon_y$ are all positive numbers to be specified.

With the sampled information $y_s(t)$, the observer now becomes 
\begin{align}\label{trobser1}
\begin{cases}
\dot{\hat{x}}\!=\!A\hat x\!+\!Bu\!+\!E_p p(\hat q\!+\!L_1(\hat y\!-\!y_s))\!+\!L_2(\hat y\!-\!y_s),\\
\hat y\!=\!C\hat x\!+\!Du,\\
\hat q\!=\!C_q\hat x,\\
\end{cases}
\end{align}
where $L_1,L_2$ are matrices to be designed.

The observer-based controller $u(t)$ has the form shown in \eqref{inputform2} where $\hat x_s(t)$ is updated at time instances $t_1^u,t_2^u,...$ by
\begin{align*}
\hat x_s(t)=\hat x(t_k^u),\;\forall t\in[t_k^u,t^u_{k+1}).
\end{align*}  
Here, $t_0^u=0$ and the triggering times $t_1^u,t_2^u,\dots$ are determined by the following triggering rule:
\begin{align}
&t_{k+1}^u=\inf\{t\mid t\geq t_k^u+\tau_u,\;\|\hat x_e(t)\|> \sigma_u \|\hat x(t)\|+\epsilon_u\}\label{triggerinput1}
\end{align}
where $\hat x_e(t)=\hat x(t_k)-\hat x(t)$ and $\tau_u,\sigma_u,\epsilon_u$ are all positive numbers to be specified. Note that the  information of $\hat x$ and $\hat x_e$  are available from the proposed observer.

The time-updating rule \eqref{triggeroutput1} (or \eqref{triggerinput1}) provides a built-in positive lower bound  $\tau_y$ (or $\tau_u$) for inter-execution times $\{t_{k+1}^y-t_{k}^y\}$ (or $\{t_{k+1}^u-t_{k}^u\}$), implying that the Zeno phenomenon will not occur. Although there is no bound guarantee on the inter-execution times between $t_{k}^y$ and $t_{k}^u$, this will not cause a problem since these two ETMs are implemented separately.

Since $y_e(t)=y_s(t)-y(t)$, the closed-loop system that combines system \eqref{dyn1}-\eqref{eq:delQC}, observer  \eqref{trobser1} and event-triggered controller  \eqref{inputform2}  is expressed compactly as
\begin{align}\label{closedtrigger3}
\dot z\!=\!A_cz\!+\!H_1p\!+\!H_2\delta  \tilde p\!+\!H_3\delta\hat p\!+\!H_4w\!+\!H_5\hat x_e\!+\!H_6y_e
\end{align}
where  $A_c$, $H_1$, $H_2$, $H_3$, $H_4$, $H_5$ are given in \eqref{eqA}, \eqref{H1},  \eqref{H5}, respectively, $\delta\hat p$ is given in \eqref{deltahatp}, and
\begin{align*}
\delta \tilde p&=p(q+\delta \tilde q)-p(q),\\
\delta \tilde q&=-(C_q+L_1C)e-L_1F_ww-L_1y_e,\\
H_6&=
\begin{pmatrix}
{\bf 0}\\
L_2
\end{pmatrix}.    
\end{align*}

\begin{theorem}\label{thm4}
	Consider the configuration shown in Figure \ref{figConfig2}  where the plant is described by \eqref{dyn1}-\eqref{eq:delQC} with $D={\bf 0}$, $p({\bf 0})={\bf 0}$, and $\|w\|_\infty\leq \omega_0$ with $\omega_0$ a positive number. Suppose that there exists $\ell>0$ such that $\|p(r)-p(s)\|\leq \ell\|r-s\|$ for any $r,s$. Suppose that there exist constants $\alpha_0>0,\mu>0$, and matrices $P\succ 0,K_1,K_2,L_1,L_2$ such that the closed-loop system \eqref{eqthm3-1}  with  controller \eqref{input1} and  observer  \eqref{obser1} satisfies  $\dot V\leq -\alpha_0 V+ \mu\|w\|^2$ where $V=z^\top Pz$. Choose any $\epsilon_y,\epsilon_u>0$, and 
		\begin{align}
		\sigma_y&=\frac{a_2\varrho \alpha_0\lambda_{m}(P)}{2\|C\|s_2},\;\;
		\sigma_u=\frac{a_1\varrho \alpha_0\lambda_{m}(P)}{2\sqrt{2}s_1},\label{eqsigmau}
		\end{align}
		where $0<\varrho<1$,
		$s_1=\|PH_5\|+\ell\|PH_3\|\|C_q\|$, $s_2=\|PH_6\|+\ell\|PH_2\|\|L_1\|$, and $a_1,a_2$ are two constants satisfying $0< a_1,a_2< 1$ and $a_1+a_2=1$. Choose $\tau_u>0$ as the solution to  the equation $\phi_1(\tau_u)=1$ where $\phi_1$ is the solution of the following ODE 
		$$\dot \phi_1=\sqrt{2}(1+\sigma_u\phi_1)(\eta_5+\eta_2\phi_1+d_1\eta_7),\;\phi_1(0)=0
		$$ 
		and choose $\tau_y>0$ as the solution to  the equation $\phi_2(\tau_y)=1$ where $\phi_2$ is the solution of the following ODE $$
		\dot \phi_2=\|C\|(1+\sigma_y\phi_2)(\eta_6+\eta_7\phi_2+d_2\eta_2),\;\phi_2(0)=0
		$$ 
		where 
		\begin{align}
		\begin{cases}
		\eta_5=\frac{\eta_1}{\sqrt{2}\sigma_u}+\frac{\eta_2\omega_0}{\epsilon_u},\\
		\eta_6=\frac{\eta_1}{\sigma_y\|C\|}+\frac{\eta_3\omega_0}{\epsilon_y},\\
		\eta_7=\ell\|H_2\|\|L_1\|+\|H_6\|,\\ d_1=\max\{\frac{\epsilon_y}{\epsilon_u},\frac{\sigma_y\|C\|}{\sqrt{2}\sigma_u}\},\\
		d_2=\max\{\frac{\epsilon_u}{\epsilon_y},\frac{\sqrt{2}\sigma_u}{\sigma_y\|C\|}\},
		\end{cases}
		\end{align}
		and $\eta_1,\eta_2,\eta_3$ are given in \eqref{eta1}. 
		Then, the closed-loop system \eqref{closedtrigger3} that impllements triggering rules \eqref{triggeroutput1} and \eqref{triggerinput1} is ISpS w.r.t. $w$.
\end{theorem}
\begin{proof}
	If the derivative of $V$ along the trajectory of \eqref{eqthm3-1} satisfies $\dot V\leq -\alpha_0 V+ \mu\|w\|^2$, then the derivative of $V$ along the trajectory of the closed-loop system \eqref{closedtrigger3} satisfies
	\begin{align}
	\dot V&\leq -\alpha_0 V+ \mu\|w\|^2+2z^\top P\Big[H_2(\delta \tilde p-\delta p)\nonumber\\
	&\qquad\qquad+H_3(\delta\hat p-\Delta p)+H_5\hat x_e+H_6y_e\Big]\nonumber\\
	&\leq -(1-\varrho)\alpha_0\lambda_{m}(P)\|z\|^2+ \mu\|w\|^2+\|z\|\Big(2s_1\|\hat x_e\|\nonumber\\
	&\; +2s_2\|y_e\|-\varrho\alpha_0\lambda_{m}(P)\|z\|\Big)\label{dotVthm4}
	\end{align}
	where the following facts are used:  
	$\|\delta \tilde p-\delta p\|\leq \ell\|\delta\tilde q-\delta q\|\leq \ell \|L_1\|\|y_e\|$,  
	$\|\delta \hat p-\Delta p\|\leq \ell\|C_q\|\|\hat x_e\|$. 
	As $\|z\|\geq \|\hat x\|/\sqrt{2}$ and $\|z\|\geq \|x\|\geq \|y\|/\|C\|$, we have 
	\begin{align}
	\|z\|\geq \frac{a_1\|\hat x\|}{\sqrt{2}}+\frac{a_2\|y\|}{\|C\|}.\label{znorm}
	\end{align}
	
	From \eqref{dotVthm4} and \eqref{znorm} we have 
	\begin{align*}
	&\dot V\leq -(1-\varrho)\alpha_0\lambda_{m}(P)\|z\|^2 + \mu\|w\|^2+ \|z\|\Big[2s_1\|\hat x_e\|\\
	&-\frac{a_1\varrho \alpha_0\lambda_{m}(P)}{\sqrt{2}}\|\hat x\|\Big]+ \|z\|\Big[2s_2\|y_e\|-\frac{a_2\varrho \alpha_0\lambda_{m}(P)}{\|C\|}\|y\|\Big].
	\end{align*}

	The condition $\|\hat x_e\|\leq \sigma_u \|\hat x\|+\epsilon_u$ 
	implies 
	\begin{align}
	\|\hat x_e\|\leq \sqrt{2}\sigma_u \|z\|+\epsilon_u,\label{triggerconu2}
	\end{align}
	and the condition $\|y_e\|\leq \sigma_y \|y\|+\epsilon_y$ 
	implies 
	\begin{align}
	\|y_e\|\leq \sigma_y\|C\| \|z\|+\epsilon_y.\label{triggercony2}
	\end{align}
	
	As long as \eqref{triggerconu2} and \eqref{triggercony2} hold, we have
	\begin{align}
	\dot V
	&\leq -[(1-\varrho)\alpha_0\lambda_{m}(P)-c]\|z\|^2+\mu\|w\|^2+ \frac{\epsilon_0^2}{4c}\label{triggerLayu2}
	\end{align}
	where $\epsilon_0=2(s_1\epsilon_u+s_2\epsilon_y)$, and $c$ is a constant satisfying $0<c<(1-\varrho)\alpha_0\lambda_{m}(P)$. 
	
	Since $\|p\|\leq \ell\|C_q\|\|x\|$, $\|\delta \hat p\|\leq \ell \|C_q\|(\|\hat x_e\|+\|e\|)$, and $\|\delta \tilde p\|\leq \ell(\|C_q+L_1C\|\|e\|+\|L_1F_w\|\|w\|+\|L_1\|\|y_e\|)$,  from \eqref{closedtrigger3} we have $\|\dot z\|\leq \eta_1\|z\|+\eta_2\|\hat x_e\|+\eta_3\|w\|+\eta_7\|y_e\|$.

	Similar to the argument in the proof of Theorem \ref{thm3}, we can show the following inequality holds when $\|\hat x_e\|\neq 0$ and $\|z\|\neq 0$:
	\begin{align}
	&\frac{\diff }{\diff t}(\frac{\|\hat x_e\|}{\sqrt{2}\sigma_u\|z\|+\epsilon_u})\leq \sqrt{2}(1+\frac{\sigma_u\|\hat x_e\|}{\sqrt{2}\sigma_u\|z\|+\epsilon_u})\times \nonumber\\
	& (\eta_5+\frac{\eta_2\|\hat x_e\|}{\sqrt{2}\sigma_u\|z\|+\epsilon_u}+\frac{\eta_7\|y_e\|}{\sqrt{2}\sigma_u\|z\|+\epsilon_u}).
	\label{dotode2}
	\end{align}
	It is easy to verify that
	$
	\frac{\eta_7\|y_e\|}{\sqrt{2}\sigma_u\|z\|+\epsilon_u}\leq d_1\frac{\eta_7\|y_e\|}{\sigma_y\|C\|\|z\|+\epsilon_y}.
	$
	Hence, from \eqref{dotode2} we have
	\begin{align*}
	&\frac{\diff }{\diff t}(\frac{\|\hat x_e\|}{\sqrt{2}\sigma_u\|z\|+\epsilon_u})\leq \sqrt{2}(1+\frac{\sigma_u\|\hat x_e\|}{\sqrt{2}\sigma_u\|z\|+\epsilon_u})\times \nonumber\\
	& (\eta_5+\frac{\eta_2\|\hat x_e\|}{\sqrt{2}\sigma_u\|z\|+\epsilon_u}+d_1\frac{\eta_7\|y_e\|}{\sigma_y\|C\|\|z\|+\epsilon_y}).
	\end{align*}
	When $\|\hat x_e\|= 0$ or $\|z\|= 0$, the upper right-hand
	derivative of $\frac{\|\hat x_e\|}{\sqrt{2}\sigma_u\|z\|+\epsilon_u}$ can be calculated similar to the proof of Theorem \ref{thm3}, which can still be captured by the inequality above.
	
	Since $\|\dot y_e\|=\|\dot y\|\leq \|C\|\|\dot x\|\leq \|C\|\|\dot z\|$, we can show that the following inequality holds  using arguments similar to those used above:
	\begin{align*}
	&\frac{\diff }{\diff t}(\frac{\|y_e\|}{\sigma_y\|C\|\|z\|+\epsilon_y})\leq \|C\|(1+\frac{\sigma_y\|y_e\|}{\sigma_y\|C\|\|z\|+\epsilon_y})\times \nonumber\\
	& (\eta_6+\frac{\eta_7\|y_e\|}{\sigma_y\|C\|\|z\|+\epsilon_y}+d_2\frac{\eta_2\|\hat x_e\|}{\sqrt{2}\sigma_u\|z\|+\epsilon_u})
	\end{align*}
	where the discussion on using the upper right-hand
	derivative is omitted since it  is similar to that used in the proof of Theorem \ref{thm3}.

	It is not hard to show that the time it takes for $\|\hat x_e\|$ (resp. $\|y_e\|$) to evolve from $0$ to $\sqrt{2}\sigma_u\|z\|+\epsilon_u$ (resp. $\sigma_y\|C\|\|z\|+\epsilon_y$) is lower bounded by $\tau_u$ (resp. $\tau_y$), which implies that \eqref{triggerconu2} holds during $[t_k^u,t_k^u+\tau_u)$, and  \eqref{triggercony2} holds during $[t_k^y,t_k^y+\tau_y)$, for any $k\geq 0$. Recalling that $\|\hat x_e\|\leq \sigma_u \|\hat x\|+\epsilon_u$ implies \eqref{triggerconu2} and $\|y_e\|\leq \sigma_y \|y\|+\epsilon_y$ implies \eqref{triggercony2}, the triggering rules \eqref{triggeroutput1} and \eqref{triggerinput1} guarantee that \eqref{triggerconu2} holds during the interval $[t_k^u, t_{k+1}^u)$ for any $k\geq 0$, and \eqref{triggercony2} holds during the interval $[t_k^y, t_{k+1}^y)$ for any $k\geq 0$. Hence, \eqref{triggerLayu2} holds for any $t\geq 0$, implying that the function $V$ is an ISpS-Lyapunov function since it satisfies \eqref{ISSineq} with $\gamma(\|z\|)=[(1-\varrho)\alpha_0\lambda_{m}(P)-c]\|z\|^2\in\mathcal{K}_\infty$, $\chi(\|w\|)=\mu\|w\|^2\in\mathcal{K}$ and $d=\epsilon_0^2/4c>0$. The conclusion follows by Proposition \ref{proISS}.
\end{proof}

\begin{remark}
The assumption $\dot V\leq -\alpha_0 V+ \mu\|w\|^2$ in Theorem \ref{thm3} and \ref{thm4} can be verified by using Theorem \ref{thm1}.
Therefore, Theorem \ref{thm1}, \ref{thm3} and \ref{thm4} altogether 
provide a systematic and constructive approach to design observer-based event-triggered controllers. One limitation of this ETC design, however, is that it  relies on the global Lipschitz constant which is normally very conservatively computed.  
\end{remark}

	\begin{remark}
		Similar to Remark \ref{remarktradeoff}, there are trade-offs in choosing parameters in triggering rules \eqref{triggeroutput1} and \eqref{triggerinput1}; for example, smaller $\epsilon_u,\epsilon_y$ reduces the ultimate bounds but decreases the inter-execution times.
	\end{remark}

\section{Simulation Example}\label{sec:example} 

\begin{figure}[!hbt]
	\centering
	\includegraphics[width=0.8\columnwidth]{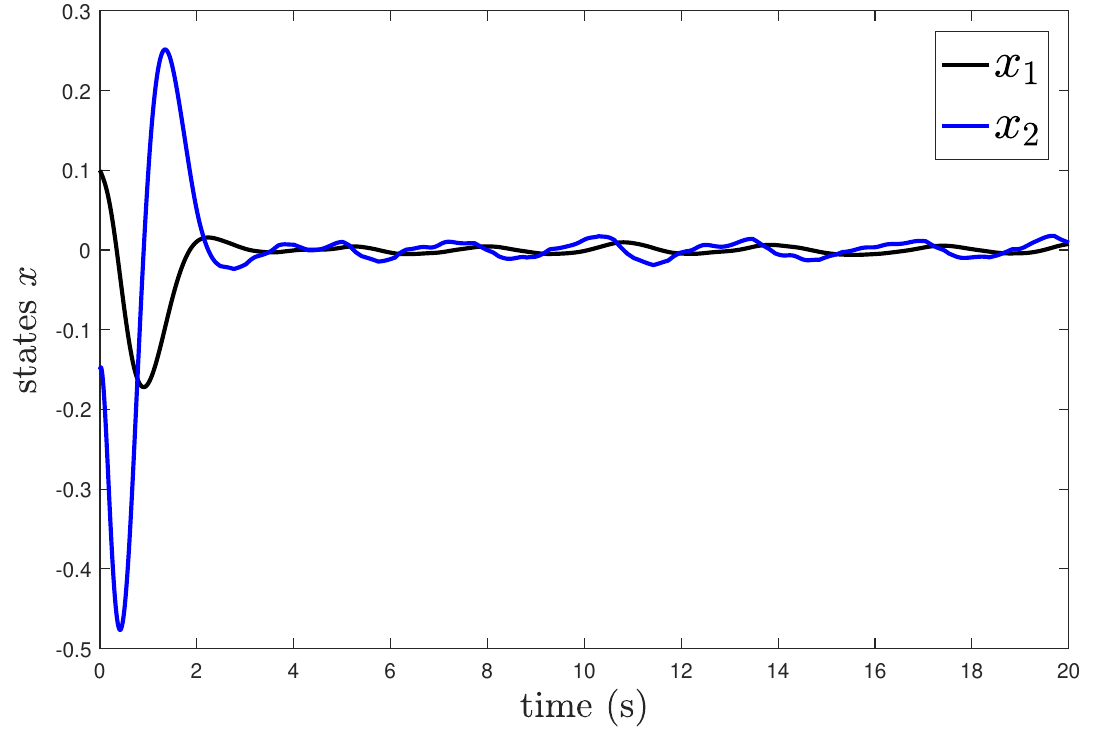}
	\caption{Trajectory of the plant state $x$.
	}\label{figstate}
\end{figure}

\begin{figure}[!hbt]
	\centering
	\includegraphics[width=0.8\columnwidth]{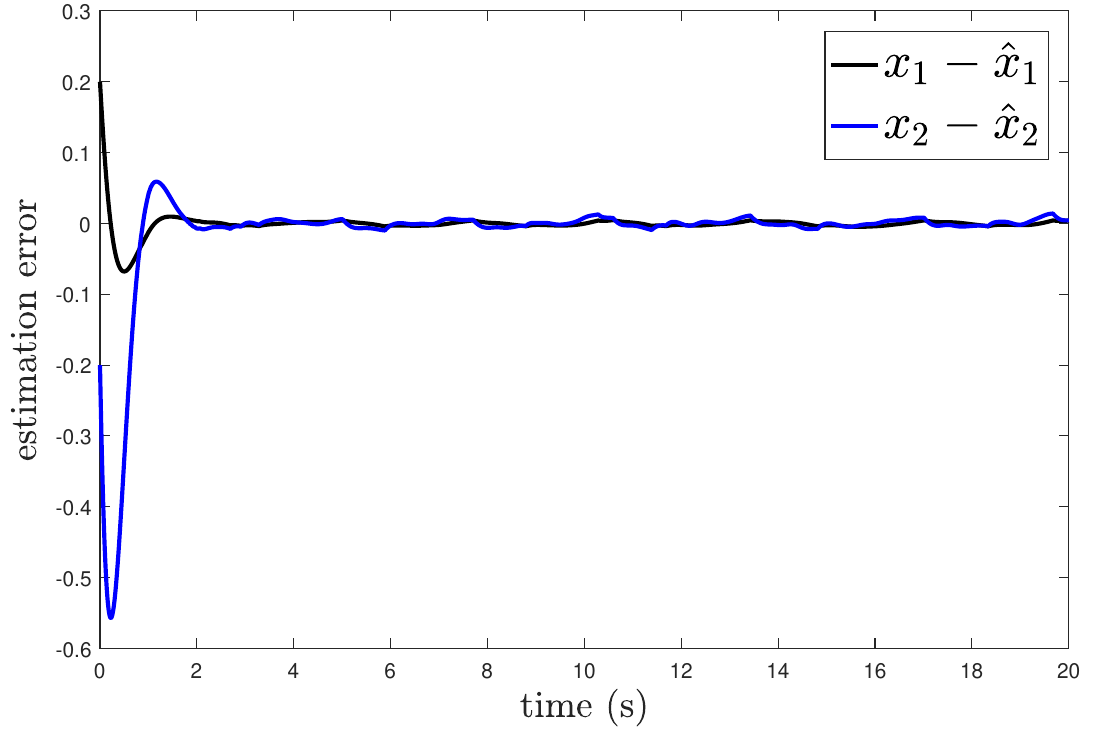}
	\caption{Trajectory of the estimation error $e=x-\hat x$.
	}\label{figerror}
\end{figure}

\begin{figure}[!hbt]
	\centering
	\includegraphics[width=0.8\columnwidth]{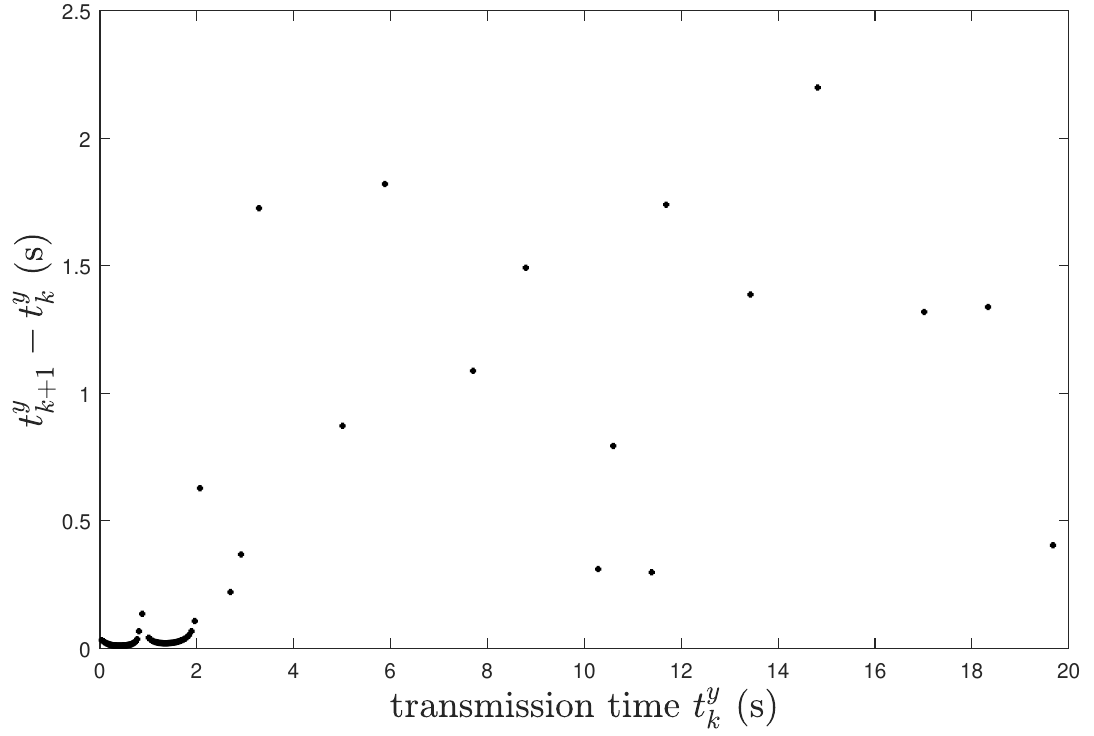}
	\caption{Inter-execution times $\{t_{k+1}^y-t_{k}^y\}$ in the observer ETM \eqref{triggeroutput1}. }\label{figtimehaty}
\end{figure}

\begin{figure}[!hbt]
	\centering
	\includegraphics[width=0.8\columnwidth]{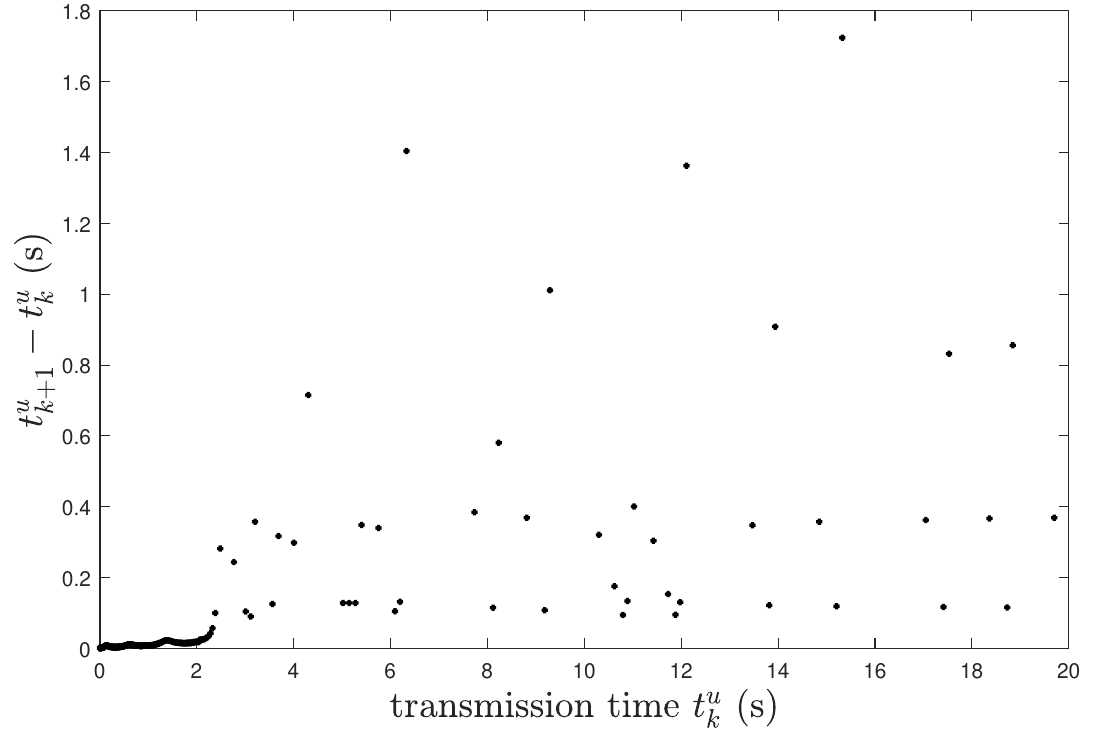}
	\caption{Inter-execution times  $\{t_{k+1}^u-t_{k}^u\}$ in the controller ETM \eqref{triggerinput1}.}\label{figtimehatx}
\end{figure}

\begin{figure}[!hbt]
	\centering
	\includegraphics[width=0.8\columnwidth]{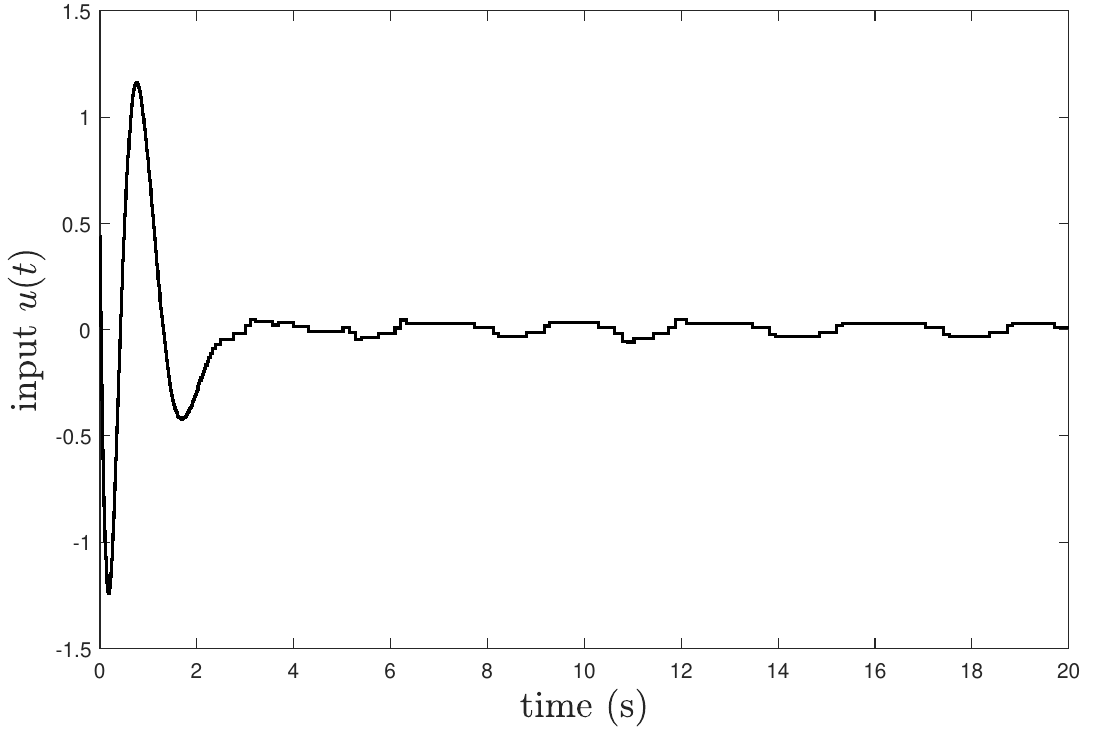}
	\caption{Trajectory of the input $u(t)$.}\label{figinput}
\end{figure}

In this section, we use a single-link robot arm example given in \cite{abdelrahim2017robust} and the configuration of Figure \ref{figConfig2}  to illustrate Theorem \ref{thm4}. 
Dynamics of the single-link robot arm are expressed as:
\begin{align*}
\dot x_1&=x_2,\\
\dot x_2&=-\sin(x_1)+u+w,\\
y&=x_1,
\end{align*}
where $x=(x_1,x_2)^\top$ is the state representing the angle and the rotational velocity, $u$ is the input representing the torque, and $w$ is the external disturbance. The system can be written in the form of \eqref{dyn1} with $A=
\begin{pmatrix}
0&1\\
0&0
\end{pmatrix}$, $B=
\begin{pmatrix}
0\\
1
\end{pmatrix}$, $C=(1,0)$, $D=0$, $E=
\begin{pmatrix}
0\\
-1
\end{pmatrix}$, $E_w=
\begin{pmatrix}
0\\
1
\end{pmatrix}$, $F_w=0$, $C_q=(1,0)$ and $p(q)=\sin(q)$. The nonlinearity $p$ satisfies  \eqref{eq:delQC} with $M=
\begin{pmatrix}
1&0\\
0&-1
\end{pmatrix}$. Recalling Remark \ref{remark2}, $p$ satisfies Assumption \ref{ass1} and \ref{ass2} with $T_1=T_2=
\begin{pmatrix}
I&{\bf 0}\\
{\bf 0}&I
\end{pmatrix}$ and 
$\mathcal{N}_1=\{(\lambda_1 I,\lambda_1 I)|\lambda_1>0\}$, $\mathcal{N}_2=\{(\frac{1}{\lambda_2} I,\frac{1}{\lambda_2} I)|\lambda_2>0\}$, which means that $X_1=Y_1=\lambda_1I$, $X_2=Y_2=\lambda_2I$. 
Additionally, the corresponding $M_{24}=-1<0$. By letting $\alpha_1=\alpha_2=1$, $\mu_1=\mu_2=0.1$, the LMIs \eqref{LMI1}-\eqref{LMI2} with variables $\lambda_1,\lambda_2,P_1,P_2$ are feasible, from which we can obtain matrix gains $L_1=-1$, $L_2=
\begin{pmatrix}
-5.1294\\
-18.0352
\end{pmatrix}$, $K_1=(-7.3936,-3.9937)$, $K_2=1$. The observer is given in \eqref{trobser1} with $L_1,L_2$ above, and the controller is given in \eqref{inputform2} with $K_1,K_2$ above. We then let $\alpha_0=0.25$, $w_0=0.02$ and recompute $P$ via \eqref{recomLMI} with the objective to be minimizing the condition number of $P$.  
With $\varrho=0.8$, $a_1=a_2=0.5$, $\epsilon_u=\epsilon_y=0.005$, we can calculate that  $\sigma_y=0.0017$, $\sigma_u=0.0023$, and $\tau_u\geq 1.07\times 10^{-4}$ s, $\tau_u\geq 7.68\times 10^{-5}$ s. In the simulations, we suppose that the random disturbance $w$ is uniformly generated from $[-w_0,w_0]$, and the initial conditions of the plant and the observer are $(0.1, -0.15)$ and $(-0.1, 0.05)$, respectively. The simulation results are shown in Figure \ref{figstate} through Figure \ref{figinput}.

Figure \ref{figstate} and Figure \ref{figerror} show trajectories of the state $x$ and the estimation error $e$, respectively. Both $x$ and $e$ eventually enter a small neighborhood of the origin as expected. 
Figure \ref{figtimehaty} shows inter-execution times  $\{t_{k+1}^y-t_{k}^y\}$ in the observer ETM \eqref{triggeroutput1}, and  Figure \ref{figtimehatx} shows inter-execution times $\{t_{k+1}^u-t_{k}^u\}$ in the controller ETM \eqref{triggerinput1}. Figure \ref{figinput} shows the trajectory of the piecewise constant input $u(t)$ that is fed into the plant. It is readily seen that the control input $u(t)$ updates its values at each sampling time $t=t_k^u$, which is determined by the triggering rule \eqref{triggerinput1}.

Denote $\tau^{min}_{[T_1,T_2]}$ and $\tau^{avg}_{[T_1,T_2]}$ as the minimal and average inter-execution times during the time interval $[T_1,T_2]$, respectively. The values of $\tau_{min}^{[0,20]}$, $\tau_{avg}^{[0,20]}$, $\tau_{min}^{[3,20]}$, $\tau_{avg}^{[3,20]}$  for the observer ETM and the controller ETM are summarized in Table \ref{tab1}. We notice that after 3 seconds, the controller input is updated about every 0.36 seconds on average, and the plant output is updated about every 1.09 seconds on average, which shows the effectiveness of our control design.

{\renewcommand{\arraystretch}{1.5}%
	\begin{table}[!hbt]\caption{Minimal and average inter-execution times for  observer  and controller ETMs }\label{tab1}
		\centering
		\begin{tabular}{c|c|c|c|c}
			\hline
			& $\tau^{min}_{[0,20]}$ & $\tau^{avg}_{[0,20]}$ & $\tau^{min}_{[3,20]}$ & $\tau^{avg}_{[3,20]}$ \\ 
			\hline
			{\footnotesize Observer ETM} & 0.0106 s & 0.1945 s& 0.2104 s& 1.0977 s\\ 
			\hline
			{\footnotesize Controller ETM} & 0.0013 s & 0.0663 s & 0.0903 s& 0.3665 s\\
			\hline
		\end{tabular} 
	\end{table}
}

\section{Conclusion}\label{sec:conclusion}

In this paper, we studied observer-based, global stabilizing control design for  incrementally quadratic nonlinear systems affected by external disturbances and measurement noise. 
We proposed LMI-based sufficient conditions for the simultaneous design of the observer and the controller in the continuous-time domain for two parameterizations of the  incremental multiplier matrices. Based on that, we investigated ETM design within the observer-based controller setting for globally Lipschitz  systems.  
The simulation  example showed the effectiveness of the controller design and the proposed triggering rule. 

\bibliographystyle{IEEEtran}
\bibliography{./TAC19-0658Bib}

\begin{thebibliography}{10}
\providecommand{\url}[1]{#1}
\csname url@samestyle\endcsname
\providecommand{\newblock}{\relax}
\providecommand{\bibinfo}[2]{#2}
\providecommand{\BIBentrySTDinterwordspacing}{\spaceskip=0pt\relax}
\providecommand{\BIBentryALTinterwordstretchfactor}{4}
\providecommand{\BIBentryALTinterwordspacing}{\spaceskip=\fontdimen2\font plus
\BIBentryALTinterwordstretchfactor\fontdimen3\font minus
  \fontdimen4\font\relax}
\providecommand{\BIBforeignlanguage}[2]{{%
\expandafter\ifx\csname l@#1\endcsname\relax
\typeout{** WARNING: IEEEtran.bst: No hyphenation pattern has been}%
\typeout{** loaded for the language `#1'. Using the pattern for}%
\typeout{** the default language instead.}%
\else
\language=\csname l@#1\endcsname
\fi
#2}}
\providecommand{\BIBdecl}{\relax}
\BIBdecl

\bibitem{krener1983linearization}
A.~J. Krener and A.~Isidori, ``Linearization by output injection and nonlinear
  observers,'' \emph{Systems \& Control Letters}, vol.~3, no.~1, pp. 47--52,
  1983.

\bibitem{vidyasagar1980stabilization}
M.~Vidyasagar, ``On the stabilization of nonlinear systems using state
  detection,'' \emph{IEEE Transactions on Automatic Control}, vol.~25, no.~3,
  pp. 504--509, 1980.

\bibitem{krstic1995nonlinear}
M.~Krstic, I.~Kanellakopoulos, and P.~V. Kokotovic, \emph{Nonlinear and
  adaptive control design}.\hskip 1em plus 0.5em minus 0.4em\relax Wiley, 1995.

\bibitem{andrieu2009unifying}
V.~Andrieu and L.~Praly, ``A unifying point of view on output feedback designs
  for global asymptotic stabilization,'' \emph{Automatica}, vol.~45, no.~8, pp.
  1789--1798, 2009.

\bibitem{borri2017luenberger}
A.~Borri, F.~Cacace, A.~De~Gaetano, A.~Germani, C.~Manes, P.~Palumbo,
  S.~Panunzi, and P.~Pepe, ``Luenberger-like observers for nonlinear time-delay
  systems with application to the artificial pancreas: The attainment of good
  performance,'' \emph{IEEE Control Systems Magazine}, vol.~37, no.~4, pp.
  33--49, 2017.

\bibitem{cheah2014observer}
C.~C. Cheah, X.~Li, X.~Yan, and D.~Sun, ``Observer-based optical manipulation
  of biological cells with robotic tweezers,'' \emph{IEEE Transactions on
  Robotics}, vol.~30, no.~1, pp. 68--80, 2014.

\bibitem{mendez2007selected}
H.~O. M{\'e}ndez-Acosta, R.~Femat, and V.~Gonz{\'a}lez-{\'A}lvarez,
  \emph{Selected topics in dynamics and control of chemical and biological
  processes}.\hskip 1em plus 0.5em minus 0.4em\relax Springer, 2007, vol. 361.

\bibitem{khalil2017high}
H.~K. Khalil, ``High-gain observers in feedback control: Application to
  permanent magnet synchronous motors,'' \emph{IEEE Control Systems Magazine},
  vol.~37, no.~3, pp. 25--41, 2017.

\bibitem{rajamani2017observers}
R.~Rajamani, Y.~Wang, G.~D. Nelson, R.~Madson, and A.~Zemouche, ``Observers
  with dual spatially separated sensors for enhanced estimation: Industrial,
  automotive, and biomedical applications,'' \emph{IEEE Control Systems
  Magazine}, vol.~37, no.~3, pp. 42--58, 2017.

\bibitem{su2007simple}
Y.~Su, P.~C. Muller, and C.~Zheng, ``A simple nonlinear observer for a class of
  uncertain mechanical systems,'' \emph{IEEE Transactions on Automatic
  Control}, vol.~52, no.~7, pp. 1340--1345, 2007.

\bibitem{ouassaid2012observer}
M.~Ouassaid, M.~Maaroufi, and M.~Cherkaoui, ``Observer-based nonlinear control
  of power system using sliding mode control strategy,'' \emph{Electric Power
  Systems Research}, vol.~84, no.~1, pp. 135--143, 2012.

\bibitem{mahmud2012full}
M.~A. Mahmud, H.~Pota, and M.~Hossain, ``Full-order nonlinear observer-based
  excitation controller design for interconnected power systems via exact
  linearization approach,'' \emph{International Journal of Electrical Power \&
  Energy Systems}, vol.~41, no.~1, pp. 54--62, 2012.

\bibitem{hong2008distributed}
Y.~Hong, G.~Chen, and L.~Bushnell, ``Distributed observers design for
  leader-following control of multi-agent networks,'' \emph{Automatica},
  vol.~44, no.~3, pp. 846--850, 2008.

\bibitem{ahmed2012high}
T.~Ahmed-Ali and F.~Lamnabhi-Lagarrigue, ``High gain observer design for some
  networked control systems,'' \emph{IEEE Transactions on Automatic Control},
  vol.~57, no.~4, pp. 995--1000, 2012.

\bibitem{kokotovic1992joy}
P.~V. Kokotovic, ``The joy of feedback: nonlinear and adaptive,'' \emph{IEEE
  Control Systems Magazine}, vol.~12, no.~3, pp. 7--17, 1992.

\bibitem{mazenc1994global}
F.~Mazenc, L.~Praly, and W.~Dayawansa, ``Global stabilization by output
  feedback: examples and counterexamples,'' \emph{Systems \& Control Letters},
  vol.~23, no.~2, pp. 119--125, 1994.

\bibitem{kokotovic2001constructive}
P.~Kokotovi{\'c} and M.~Arcak, ``Constructive nonlinear control: a historical
  perspective,'' \emph{Automatica}, vol.~37, no.~5, pp. 637--662, 2001.

\bibitem{esfandiari1992output}
F.~Esfandiari and H.~K. Khalil, ``Output feedback stabilization of fully
  linearizable systems,'' \emph{International Journal of control}, vol.~56,
  no.~5, pp. 1007--1037, 1992.

\bibitem{khalil2014high}
H.~K. Khalil and L.~Praly, ``High-gain observers in nonlinear feedback
  control,'' \emph{International Journal of Robust and Nonlinear Control},
  vol.~24, no.~6, pp. 993--1015, 2014.

\bibitem{teel1994global}
A.~Teel and L.~Praly, ``Global stabilizability and observability imply
  semi-global stabilizability by output feedback,'' \emph{Systems \& Control
  Letters}, vol.~22, no.~5, pp. 313--325, 1994.

\bibitem{khalil1993semiglobal}
H.~Khalil and F.~Esfandiari, ``Semiglobal stabilization of a class of nonlinear
  systems using output feedback,'' \emph{IEEE Transactions on Automatic
  Control}, vol.~38, no.~9, pp. 1412--1415, 1993.

\bibitem{atassi2000separation}
A.~Atassi and H.~Khalil, ``Separation results for the stabilization of
  nonlinear systems using different high-gain observer designs,'' \emph{Systems
  \& Control Letters}, vol.~39, no.~3, pp. 183--191, 2000.

\bibitem{atassi2001separation}
------, ``A separation principle for the control of a class of nonlinear
  systems,'' \emph{IEEE Transactions on Automatic Control}, vol.~46, no.~5, pp.
  742--746, 2001.

\bibitem{gauthier1992separation}
J.~Gauthier and I.~Kupka, ``A separation principle for bilinear systems with
  dissipative drift,'' \emph{IEEE transactions on automatic control}, vol.~37,
  no.~12, pp. 1970--1974, 1992.

\bibitem{tsinias1993sontag}
J.~Tsinias, ``Sontag's ‘input to state stability condition’and global
  stabilization using state detection,'' \emph{Systems \& Control Letters},
  vol.~20, no.~3, pp. 219--226, 1993.

\bibitem{lin1995bounded}
W.~Lin, ``Bounded smooth state feedback and a global separation principle for
  non-affine nonlinear systems,'' \emph{Systems \& Control Letters}, vol.~26,
  no.~1, pp. 41--53, 1995.

\bibitem{arcak2001nonlinear}
M.~Arcak and P.~Kokotovi{\'c}, ``Nonlinear observers: a circle criterion design
  and robustness analysis,'' \emph{Automatica}, vol.~37, no.~12, pp.
  1923--1930, 2001.

\bibitem{arcak2001observer}
M.~Arcak and P.~Kokotovic, ``Observer-based control of systems with
  slope-restricted nonlinearities,'' \emph{IEEE Transactions on Automatic
  Control}, vol.~46, no.~7, pp. 1146--1150, 2001.

\bibitem{arcak2005certainty}
M.~Arcak, ``Certainty-equivalence output-feedback design with circle-criterion
  observers,'' \emph{IEEE Transactions on Automatic Control}, vol.~50, no.~6,
  pp. 905--909, 2005.

\bibitem{fossen1999passive}
T.~I. Fossen and J.~P. Strand, ``Passive nonlinear observer design for ships
  using lyapunov methods: full-scale experiments with a supply vessel,''
  \emph{Automatica}, vol.~35, no.~1, pp. 3--16, 1999.

\bibitem{loria2000separation}
A.~Loria, T.~I. Fossen, and E.~Panteley, ``A separation principle for dynamic
  positioning of ships: Theoretical and experimental results,'' \emph{IEEE
  Transactions on Control Systems Technology}, vol.~8, no.~2, pp. 332--343,
  2000.

\bibitem{praly1993stabilization}
L.~Praly and Z.-P. Jiang, ``Stabilization by output feedback for systems with
  {ISS} inverse dynamics,'' \emph{Systems \& Control Letters}, vol.~21, no.~1,
  pp. 19--33, 1993.

\bibitem{pomet1993dynamic}
J.-B. Pomet, R.~M. Hirschorn, and W.~Cebuhar, ``Dynamic output feedback
  regulation for a class of nonlinear systems,'' \emph{Mathematics of Control,
  Signals and Systems}, vol.~6, no.~2, pp. 106--124, 1993.

\bibitem{praly1993lyapunov}
L.~Praly, ``Lyapunov design of a dynamic output feedback for systems linear in
  their unmeasured state components,'' in \emph{Nonlinear Control Systems
  Design}.\hskip 1em plus 0.5em minus 0.4em\relax Pergamon, 1993, pp. 63--68.

\bibitem{boyd1994linear}
S.~Boyd, L.~El~Ghaoui, E.~Feron, and V.~Balakrishnan, \emph{Linear matrix
  inequalities in system and control theory}.\hskip 1em plus 0.5em minus
  0.4em\relax SIAM, 1994.

\bibitem{lien2004robust}
C.-H. Lien, ``Robust observer-based control of systems with state perturbations
  via {LMI} approach,'' \emph{IEEE Transactions on Automatic Control}, vol.~49,
  no.~8, pp. 1365--1370, 2004.

\bibitem{kheloufi2013lmi}
H.~Kheloufi, A.~Zemouche, F.~Bedouhene, and M.~Boutayeb, ``On {LMI} conditions
  to design observer-based controllers for linear systems with parameter
  uncertainties,'' \emph{Automatica}, vol.~49, no.~12, pp. 3700--3704, 2013.

\bibitem{zemouche2017robust}
A.~Zemouche, R.~Rajamani, H.~Kheloufi, and F.~Bedouhene, ``Robust
  observer-based stabilization of {L}ipschitz nonlinear uncertain systems via
  {LMI}s-discussions and new design procedure,'' \emph{International Journal of
  Robust and Nonlinear Control}, vol.~27, no.~11, pp. 1915--1939, 2017.

\bibitem{toker1995np}
O.~Toker and H.~Ozbay, ``On the {NP}-hardness of solving bilinear matrix
  inequalities and simultaneous stabilization with static output feedback,'' in
  \emph{American Control Conference}.\hskip 1em plus 0.5em minus 0.4em\relax
  IEEE, 1995, pp. 2525--2526.

\bibitem{wang2018sequential}
Y.~Wang, R.~Rajamani, and A.~Zemouche, ``Sequential {LMI} approach for the
  design of a {BMI}-based robust observer state feedback controller with
  nonlinear uncertainties,'' \emph{International Journal of Robust and
  Nonlinear Control}, vol.~28, no.~4, pp. 1246--1260, 2018.

\bibitem{kim2017robust}
K.-K.~K. Kim and R.~D. Braatz, ``Robust static and fixed-order dynamic output
  feedback control of discrete-time parametric uncertain lur{\'e} systems:
  Sequential sdp relaxation approaches,'' \emph{Optimal Control Applications
  and Methods}, vol.~38, no.~1, pp. 36--58, 2017.

\bibitem{grandvallet2013new}
B.~Grandvallet, A.~Zemouche, H.~Souley-Ali, and M.~Boutayeb, ``New {LMI}
  condition for observer-based {H}$_\infty$ stabilization of a class of
  nonlinear discrete-time systems,'' \emph{SIAM Journal on Control and
  Optimization}, vol.~51, no.~1, pp. 784--800, 2013.

\bibitem{ekram2017observer}
M.~Ekramian, ``Observer-based controller for {L}ipschitz nonlinear systems,''
  \emph{International Journal of Systems Science}, vol.~48, no.~16, pp.
  3411--3418, 2017.

\bibitem{megretski1997system}
A.~Megretski and A.~Rantzer, ``System analysis via integral quadratic
  constraints,'' \emph{IEEE Transactions on Automatic Control}, vol.~42, no.~6,
  pp. 819--830, 1997.

\bibitem{acikmesethesis}
B.~A{\c{c}}{\i}kme{\c{s}}e, ``Stabilization, observation, tracking and
  disturbance rejection for uncertain/nonlinear and time-varying systems,''
  Ph.D. dissertation, Purdue University, 2002.

\bibitem{acikmese2005observers}
B.~A{\c{c}}{\i}kme{\c{s}}e and M.~Corless, ``Observers for systems with
  nonlinearities satisfying an incremental quadratic inequality,'' in
  \emph{American Control Conference}, 2005, pp. 3622--3629.

\bibitem{accikmecse2008stability}
------, ``Stability analysis with quadratic lyapunov functions: some necessary
  and sufficient multiplier conditions,'' \emph{Systems \& Control Letters},
  vol.~57, no.~1, pp. 78--94, 2008.

\bibitem{accikmecse2011observers}
------, ``Observers for systems with nonlinearities satisfying incremental
  quadratic constraints,'' \emph{Automatica}, vol.~47, no.~7, pp. 1339--1348,
  2011.

\bibitem{alto13incremental}
L.~D'Alto and M.~Corless, ``Incremental quadratic stability,'' \emph{Numerical
  Algebra, Control and Optimization}, vol.~3, no.~1, pp. 175--201, 2013.

\bibitem{ankush2017state}
A.~Chakrabarty, M.~Corless, G.~T. Buzzard, S.~H. Zak, and A.~E. Rundell,
  ``State and unknown input observers for nonlinear systems with bounded
  exogenous inputs,'' \emph{IEEE Transactions on Automatic Control}, vol.~62,
  no.~11, pp. 5497--5510, 2017.

\bibitem{chen2007robust}
M.-S. Chen and C.-C. Chen, ``Robust nonlinear observer for {L}ipschitz
  nonlinear systems subject to disturbances,'' \emph{IEEE Transactions on
  Automatic control}, vol.~52, no.~12, pp. 2365--2369, 2007.

\bibitem{arcak1999observer}
M.~Arcak and P.~Kokotovi{\'c}, ``Observer-based stabilization of systems with
  monotonic nonlinearities,'' \emph{Asian Journal of Control}, vol.~1, no.~1,
  pp. 42--48, 1999.

\bibitem{fan2003observer}
X.~Fan and M.~Arcak, ``Observer design for systems with multivariable monotone
  nonlinearities,'' \emph{Systems \& Control Letters}, vol.~50, no.~4, pp.
  319--330, 2003.

\bibitem{heemels2012introduction}
W.~Heemels, K.~H. Johansson, and P.~Tabuada, ``An introduction to
  event-triggered and self-triggered control,'' in \emph{IEEE Conference on
  Decision and Control}, 2012, pp. 3270--3285.

\bibitem{tabuada2007event}
P.~Tabuada, ``Event-triggered real-time scheduling of stabilizing control
  tasks,'' \emph{IEEE Transactions on Automatic Control}, vol.~52, no.~9, pp.
  1680--1685, 2007.

\bibitem{postoyan2015framework}
R.~Postoyan, P.~Tabuada, D.~Ne{\v{s}}i{\'c}, and A.~Anta, ``A framework for the
  event-triggered stabilization of nonlinear systems,'' \emph{IEEE Transactions
  on Automatic Control}, vol.~60, no.~4, pp. 982--996, 2015.

\bibitem{tallapragada2012event}
P.~Tallapragada and N.~Chopra, ``Event-triggered dynamic output feedback
  control for {LTI} systems,'' in \emph{IEEE Conference on Decision and
  Control}, 2012, pp. 6597--6602.

\bibitem{tarbouriech2016observer}
S.~Tarbouriech, A.~Seuret, J.~M.~G. da~Silva~Jr, and D.~Sbarbaro,
  ``Observer-based event-triggered control co-design for linear systems,''
  \emph{IET Control Theory \& Applications}, vol.~10, no.~18, pp. 2466--2473,
  2016.

\bibitem{donkers2012output}
M.~Donkers and W.~Heemels, ``Output-based event-triggered control with
  guaranteed $\mathcal{L}_\infty$-gain and improved and decentralized
  event-triggering,'' \emph{IEEE Transactions on Automatic Control}, vol.~57,
  no.~6, pp. 1362--1376, 2012.

\bibitem{borgers2014event}
D.~N. Borgers and W.~M. Heemels, ``Event-separation properties of
  event-triggered control systems,'' \emph{IEEE Transactions on Automatic
  Control}, vol.~59, no.~10, pp. 2644--2656, 2014.

\bibitem{dolk2017output}
V.~Dolk, D.~P. Borgers, and W.~Heemels, ``Output-based and decentralized
  dynamic event-triggered control with guaranteed $\mathcal{L}_p$-gain
  performance and zeno-freeness,'' \emph{IEEE Transactions on Automatic
  Control}, vol.~62, no.~1, pp. 34--49, 2017.

\bibitem{abdelrahim2017robust}
M.~Abdelrahim, R.~Postoyan, J.~Daafouz, and D.~Ne{\v{s}}i{\'c}, ``Robust
  event-triggered output feedback controllers for nonlinear systems,''
  \emph{Automatica}, vol.~75, pp. 96--108, 2017.

\bibitem{xu2018acc}
X.~Xu, B.~A{\c{c}}{\i}kme{\c{s}}e, M.~Corless, and H.~Sartipizadeh,
  ``Observer-based output feedback control design for systems with
  incrementally conic nonlinearities,'' in \emph{American Control Conference},
  2018, pp. 1364--1369.

\bibitem{jiang1996lyapunov}
Z.-P. Jiang, I.~M. Mareels, and Y.~Wang, ``A {L}yapunov formulation of the
  nonlinear small-gain theorem for interconnected {ISS} systems,''
  \emph{Automatica}, vol.~32, no.~8, pp. 1211--1215, 1996.

\bibitem{isidori2013nonlinear}
A.~Isidori, \emph{Nonlinear Control Systems-II}.\hskip 1em plus 0.5em minus
  0.4em\relax Springer Science \& Business Media, 2013.

\bibitem{sontag1995characterizations}
E.~D. Sontag and Y.~Wang, ``On characterizations of the input-to-state
  stability property,'' \emph{Systems \& Control Letters}, vol.~24, no.~5, pp.
  351--359, 1995.

\end{thebibliography}

\end{document}